\documentclass[a4paper,11pt]{article}

\title{Regularity estimate and sparse approximation of pathwise robust Duncan-Mortensen-Zakai equation}
\author{Yuhua Meng\thanks{Department of Mathematics, The University of Hong Kong, Pokfulam Road, Hong Kong SAR, P.R.China. Email: {yvette00@connect.hku.hk.} } 
\and Zhongjian Wang\thanks{Division of Mathematical Sciences, School of Physical and Mathematical Sciences, Nanyang
Technological University, 21 Nanyang Link Singapore 637371. Email: {zhongjian.wang@ntu.edu.sg.} }
\and Stephen S.T. Yau\thanks{Hetao Institute of Mathematics and Interdisciplinary Sciences (HIMIS), Shenzhen 518000, Guangdong, P. R. China, Beijing Institute of Mathematical Sciences and Applications (BIMSA), Beijing 101408, P. R. China, Department of Mathematical Sciences, Tsinghua University, Beijing 100084, P. R. China
Email: {yau@uic.edu.} }
\and Zhiwen Zhang\thanks{Corresponding author. Department of Mathematics, The University of Hong Kong, Pokfulam Road, Hong Kong SAR, P.R.China. Materials Innovation Institute for Life Sciences and Energy (MILES), HKU-SIRI,
Shenzhen, P.R. China. Email: {zhangzw@hku.hk.})}
}

\usepackage{graphicx,subfigure}
\usepackage{xcolor}
\usepackage[shortlabels]{enumitem}
\usepackage{amsmath,amsfonts,amsthm,amssymb}
\usepackage{cleveref}
\usepackage{xcolor}
\usepackage{verbatim}
\usepackage{placeins}
\usepackage{bm}
 \usepackage{booktabs} 
 \usepackage{algorithm}
 \usepackage{algpseudocode}
 \usepackage[backend=biber, style=numeric, maxnames=99]{biblatex}
\newtheorem{theorem}{Theorem}[section]
\newtheorem{assumption}[theorem]{Assumption}
\newtheorem{proposition}[theorem]{Proposition}
\newtheorem{lemma}[theorem]{Lemma}
\newtheorem{corollary}[theorem]{Corollary}

\newtheorem{definition}[theorem]{Definition}
\newtheorem{remark}[theorem]{Remark}

\algnewcommand{\Input}{\item[\textbf{Input:}]}
\algnewcommand{\Output}{\item[\textbf{Output:}]}
\begin{document}
%\title{Regularity Estimate of pathwise robust DMZ equation and its sparse approximation}
\maketitle

% REQUIRED
\begin{abstract}
In this paper, we establish an \textit{a priori} estimate for arbitrary-order derivatives of the solution to the pathwise robust Duncan-Mortensen-Zakai (DMZ) equation within the framework of weighted Sobolev spaces. The weight function, which vanishes on the physical boundary, is crucial for the \textit{a priori} estimate, but introduces a loss of regularity near the boundary. Therefore, we employ the Sobolev inequalities and their weighted analogues to sharpen the regularity bound, providing improvements in both classical Sobolev spaces and H{\"o}lder continuity estimates. The refined regularity estimate reinforces the plausibility of the quantized tensor train (QTT) method in \cite{li2022solving} and provides convergence guarantees of the method. To further enhance the capacity of the method to solve the nonlinear filtering problem in a real-time manner, we reduce the complexity of the method under the assumption of a functional polyadic state drift $f$ and observation $h$. Finally, we perform numerical simulations to reaffirm our theory. For high-dimensional cubic sensor problems, our method demonstrates superior efficiency and accuracy in comparison to the particle filter (PF) and the extended Kalman filter (EKF). Beyond this, for multi-mode problems, while the PF exhibits a lack of precision due to its stochastic nature and the EKF is constrained by its Gaussian assumption, the enhanced method provides an accurate reconstruction of the multi-mode conditional density function. 



\noindent
\end{abstract}

% REQUIRED

% 93E11: Filtering in stochastic control theory
% 93-08: Computational methods for problems pertaining to systems and control theory

% 35B45: A priori estimates in context of PDEs
% 35k20: IBVPs for second-order parabolic equations.

% 65C30: Numerical solutions to stochastic differential and integral equations
% 65F55: Numerical methods for low-rank matrix approximations; matrix compression.
% 65M06: FDMs for IVPs and IBVPs involving PDE.
% 65M12: Stability and Convergence of Numerical methods for IVPs and IBVPs involving PDEs.
% 65M15: Error bounds for initial value and initial-boundary value problems involving PDEs
\section{Introduction}
The problem of nonlinear filtering (NLF) stems from the field of signal processing and has been important in various fields of engineering. A general signal observation model with explicit time dependence is stated as follows:
\begin{align}\label{eq:signal_observ_model}
\begin{cases}
&d{x}_t={f}(x_t,t)dt+{G}(t)d{v}_t, \\
&d{y}_t={h}(x_t,t)dt+d{w}_t,
\end{cases}
\end{align}
where ${x}_t\in\mathbb{R}^d$ is the state of the system at time $t$ with $x_0$ follows some known initial probability distribution $\sigma_0$, $y_t\in\mathbb{R}^m$ is the observation at time $t$ with $y_0=0$, ${v}_t$ and ${w}_t$ are Brownian motions with variance matrices $E(d{v}_t d{v}_t^T) ={Q}(t)dt$ and $E(d{w}_td{w}_t^T)={S}(t)dt$, respectively.

The core problem in the filtering is to estimate the underlying states of the system given noisy observations. For linear systems with Gaussian assumption, the Kalman filter (KF) \cite{Kalman1960anew, KalmanBucy1961} provides an optimal solution by recursively computing the conditional mean value and covariance of the state. However, most practical systems are nonlinear, motivating extensions of the KF to address such cases. Notable generalizations include the extended Kalman filter (EKF), the unscented Kalman filter, and the ensemble Kalman filter. Another widely used approach is the particle filter (PF), which approximates the state distribution by a set of random samples; see \cite{jazwinski2013stochastic, arulampalam2002tutorial} and the references therein. The PF is applicable to fully nonlinear and non-Gaussian systems. Nevertheless, its Monte Carlo nature can make it computationally intensive and thereby unsuitable for real-time applications.

An alternative approach models the evolution of the conditional probability density functions via the Duncan-Mortensen-Zakai (DMZ) equation, a stochastic partial differential equation (SPDE) \cref{eq:DMZ_equation} \cite{duncan1967probability, mortensen1966optimal, zakai1969optimal}. By introducing an invertible transformation, the DMZ equation is reformulated as a deterministic partial differential equation (PDE) with stochastic coefficients, known as the pathwise-robust DMZ (PR-DMZ) equation \cite{yau2000real,yau2008real}; see \cref{eq:robustDMZ_v1,eq:robustDMZ_v2}. 

The well-posedness of the DMZ and PR-DMZ equations has been thoroughly investigated. Pardoux \cite{pardouxt1980stochastic} proved the existence and uniqueness of the solution for bounded drift term $f\in C^1$ and the observe term $h\in C^2$. Fleming and Mitter \cite{fleming1982optimal} considered the case where $f$ and $\nabla f$ are bounded and $h$ is of polynomial growth. Baras, Blankenship, and Hopkins \cite{baras2003existence} proved the results for unbounded coefficients, but only in one dimension. In \cite{yau2000real}, Yau and Yau establish well-posedness with $f$ and $h$ having at most linear growth. In the appendices of \cite{yau2008real}, they further proved the existence and uniqueness in weighted Sobolev space with time-independent $f$ and $h$ under some mild growth conditions. Luo and Yau in \cite{luo2013complete} generalized the results in \cite{yau2008real} to the time-dependent case, which is stated in \cref{prop:low_priori_est}. Furthermore, the validity of constraining the PR-DMZ equation to a bounded domain with zero-Dirichlet condition has been established in \cite{yau2008real, sun2024convergence}.

Several techniques have been proposed to numerically solve the DMZ or PR-DMZ equation. The splitting-up method \cite{bensoussan1990approximation, gyongy2003splitting} decomposes the time integration into a stochastic and deterministic component. The $S^3$ algorithm \cite{lototsky1997nonlinear} develop the approach based on the Cameron-Martin version of Wiener chaos expansion. Both with convergence guaranteed only when the drift $f$, observation $h$, and diffusion term $G$ are bounded. In contrast, the Yau-Yau's algorithm \cite{yau2000real} is applicable to cases in which $f$ and $h$ under mild growth conditions, with its explicit convergence rate under a spectral method proved in \cite{Kang2025explicit}. Here, we instead develop the convergence analysis for the tensor train approximation within the finite-difference framework.   

Numerous approaches have been proposed to efficiently implement Yau-Yau's algorithm for high-dimensional problems. In \cite{luo2013HChermite}, Luo and Yau developed the Hermite spectral method (HSM) with hyperbolic cross approximation. The spectral method tackle the nonlinear filtering problem in medium-high dimensional cases, while the full scale resolution makes solving higher dimensional problems a bottleneck. Another type of algorithm exploit the sparse structures of the solution. Wang \textit{et al.} employed proper orthogonal decomposition (POD) for real-time observation synchronization \cite{wang2019proper}. Finite difference methods (FDM), combined with proper techniques, can also mitigate the curse of dimensionality. Yueh \textit{et al.} proposed a discrete sine transform-based quasi-implicit Euler method in \cite{yueh2014efficient}, solving NLF problems up to 6D. Tensor train (TT)/ Quantized tensor train (QTT) decompositions \cite{oseledets2011tensor,oseledets2010approximation,dolgov2012fast} enable sparse approximations, significantly reducing computation costs. Notably, Li \textit{et al.} demonstrated the effectiveness of QTT in solving 6D problems by the explicit Euler method \cite{li2022solving}. 

While previous research has made significant progress in attending to non-linear observations, existing regularity results \emph{up to first-order derivatives} (e.g. Theorem A.1 in \cite{luo2013complete}) are not sufficient for two key applications: convergence analysis of numerical schemes and the theoretical justification of sparse approximations to the solutions, which often require higher-order derivatives or the H{\" o}lder continuity (\cite{luo2013HChermite,bigoni2016spectral, griebel2023low}). To close this gap, we establish \textit{a priori} regularity estimates for derivatives of arbitrary order for PR-DMZ solution restricted on a bounded domain with a Dirichlet boundary condition. %Our estimate can be further sharpened via the embedding operators for the involved function spaces. 

The regularity estimate for the PR-DMZ solution derived in \cite{luo2013complete} is based on controlling the evolution of the energy function. A key challenge of this approach arises when estimating higher-order derivatives: integration by parts generates boundary terms that cannot be controlled by the Dirichlet condition. While both the \textit{a priori} estimates of the solution and its first-order derivatives benefit from the natural cancellation of such terms under the zero Dirichlet condition, higher-order estimates lack such auxiliary boundary conditions, leaving such terms a priori unbounded.

To address this, we introduce the weight functions that vanish on the boundary and satisfy appropriate growth conditions. These weights effectively remove problematic boundary contributions, enabling us to complete the energy estimate in \cref{thm:prioriesti} for derivatives of any order. The existence of such weight functions can be guaranteed in some specific cases. Specifically, the power-type weight function in \cref{eq:power_type_weight} is a typical example of the zero-trace weight functions. Under this weight functions, the assumptions in \cref{thm:prioriesti} reduce to requiring $F$ and $J$ to belong to classical Sobolev space for regular bounded domains such as spheres and cubes.

The weighted approach provides a higher-order regularity estimate, but mainly focuses on interior behaviour, and does not explicitly characterize derivative behaviour near the boundary. To bridge this gap, we apply embedding theorems for weighted Sobolev space \cite{kufner1980weighted} to improve the estimate in the vicinity of boundary. Additionally, the classical results from trace and Sobolev embedding theorems \cite{adams2003sobolev} ensure that sufficiently strong Sobolev regularity implies boundary regularity and H{\" o}lder continuity of the solution.  

On the computational side, we focus on the enhancement of QTT method \cite{li2022solving}. Although the low-rank structure of the operators in \cref{eq:FKEunbd} enables a fast construction of low-rank QTT structure (see \cite{kazeev2012low} and \cref{alg:offline_stage}), efficiently evaluating the new observation assimilation tensor in online procedure remains a challenge. The cost of converting a full array ${\bf A}\in \mathbb{R}^{n_1\times\cdots\times n_d}$ into tensor train format is $\mathcal{O}(n^dr^2)$ \cite{oseledets2010approximation}, where $n_i\le n$ and $r$ is the bound of compression rank. This exponential scaling of dimension $d$ deteriorates the efficiency and even the feasibility of QTT method. 

To overcome this limitation, fast evaluation algorithms such as those in \cite{savostyanov2011fast, savostyanov2014quasioptimality} can be applied. In this paper, we present a refined implementation of QTT-format tensor construction, where the functions $f$ and $h$ are represented in a functional-polyadic structure. This structures enables the TT/QTT construction to be carried out using only operations: Kronecker product and addition of a low-dimensional tensor. As a result, only full arrays of small size need to be converted to the TT/QTT format, leading to a substantial reduction in both storage requirements and computational complexity. Numerical experiments in cubic sensor examples demonstrate that the improved method achieves much superior computational performance while maintaining accuracy for high dimensions. Moreover, our method outperforms the particle filter (PF) and extended Kalman filter (EKF) in recovering the multi-mode (and hence non-log-concave) conditional density (see \cref{fig:bimodal}).

The remainder of this article is organized as follows. In \cref{sec:preliminaries}, we introduce the basic results on model and algorithm, the basic idea of weighted Sobolev space, and TT/QTT format. In \cref{sec:a priori estimate}, the a priori estimate is discussed in the framework of weighted Sobolev space. In \cref{sec:sparse approx}, we establish an efficient implementation of the QTT method in \cite{li2022solving} and its convergence analysis. The numerical results are presented in \cref{sec:experiments}. \Cref{sec:conclusions} concludes the article.
\medskip

\section{Preliminaries}
\label{sec:preliminaries}
In this section, we first present some preliminaries on PR-DMZ equation and Yau-Yau's algorithm. In order to solve high-dimensional PR-DMZ in a real-time manner, sparse approximation is required to save the computation. Further, the regularity of the solution is the basis of the existence of a low-rank approximation. Therefore, we then introduce the weighted Sobolev space framework \cite{kufner1980weighted} for \textit{a priori} estimate, and the tensor train format \cite{oseledets2011tensor, bigoni2016spectral} in sparse numerical solutions for completeness. 
\subsection{Pathwise Robust DMZ model}
Consider the signal processing model \cref{eq:signal_observ_model}, the unnormalized density function $\sigma(x,t)$ of $x_t$ conditioned on $\{{y}_{t'}:0\le t' \le t\}$ satisfies the DMZ equation \cref{eq:DMZ_equation}:
\begin{align}\label{eq:DMZ_equation}
\begin{cases}
&d\sigma({x},t)=\mathcal{L}\sigma({x},t) dt+\sigma({x},t){h}^T {S}^{-1} d{y}_t, \\
&\sigma({x},0)=\sigma_0(x),
\end{cases}
\end{align}
where
$$\mathcal{L}(\cdot)=\frac{1}{2}\sum_{ij=1}^d\frac{\partial^2}{\partial x_i\partial x_j}\big(({GQG}^T)_{ij}\cdot\big)-\sum_{i=1}^d \frac {\partial(f_i\cdot)} {\partial x_i}.$$

Equation \cref{eq:DMZ_equation} is an Ito stochastic partial differential equation with $dy_t$. In real applications, one observes data from only one realization of $x_t$; thus, a robust estimator that is not sensitive to observed paths is of more interest. By applying an invertible transformation to $\sigma(x,t)$ with each given observation path $y_t$

$$\sigma({x},t)=\exp({h}^T{S}^{-1}{y}_t)u({x},t),$$
the DMZ equation is reformulated into a deterministic PDE, which is referred as the pathwise-robust DMZ equation, i.e., PR-DMZ equation,
\begin{align}\label{eq:robustDMZ_v1}
    \begin{cases}
        & u_t+\frac{\partial}{\partial t}(h^TS^{-1})y_{t}u(x,t) \\&\qquad\qquad= \exp{(-h^TS^{-1}y_t)}(\mathcal{L}-\frac{1}{2}h^TS^{-1}h)(\exp{(h^TS^{-1}y_t)}u(x,t)),\\
        &u(x,0)=\sigma_0(x).
    \end{cases}
\end{align}
Or equivalently, we denote,
\begin{align}\label{eq:robustDMZ_v2}
    \begin{cases}
        &u_t=\frac{1}{2}\widetilde{\Delta}u({x},t)+F({x},t)\cdot\nabla u({x},t)+J({x},t) u({x},t),\\
        &u({x},0)=\sigma_0({x}),
    \end{cases}
\end{align}
where,
\begin{align*}
   \widetilde{\Delta}(\cdot) =&\sum_{ij}({GQG}^T)_{ij}(\frac{\partial ^2}{\partial x_i\partial x_j}\cdot),\\
    F({x},t)=&\text{div}({GQG}^T)+\widetilde{\nabla}(K)-{f},\\
    J( {x},t)=&-\frac{\partial}{\partial t}({h}^T{S}^{-1}){y}_t+\frac{1}{2}\sum_{ij}\frac{\partial^2({GQG}^T)_{ij}}{\partial {x}_i\partial x_j}+\text{div} ({GQG}^T)\cdot\nabla K+\notag
    \\&\frac{1}{2}\widetilde{\Delta}(K) + \frac{1}{2}\widetilde{\nabla}(K)\cdot\nabla K-\text{div}f-{f}\cdot\nabla K-\frac{1}{2}{h}^T{S}^{-1}{h},\\
    \widetilde{\nabla}(\cdot)=&({GQG}^T)\nabla(\cdot),\\
    K({x},t)=&{h}^T{S}^{-1}{y}_t.
\end{align*}
The \textit{a priori} estimate of the PR-DMZ solution and its first-order derivative is stated below.
\begin{proposition}[Theorem A.1 in \cite{luo2013complete}]\label{prop:low_priori_est}
     Let $u_R$ be the solution to the PR-DMZ equation \cref{eq:IBVP} on $[0,T]\times B_R$, where $B_R = \{ x\in\mathbb{R}^d:|x|\le R \}$ is a ball of radius $R$. If there exists positive functions $g_0(x),g_1(x)$ on $\mathbb{R}^d$, such that $g_0(x),g_1(x)$ and the equation satisfy some mild conditions (essentially assume that $J$ and the growth of $GQG^T$, $F$ and $g_0,\,g_1$ is upper bounded), then for $0\le t\le T$, there exists a positive constant $C$, such that
    \begin{subequations}\label{eq:low_priori_est}
    \begin{align}
        &\int_{B_R}e^{2g_0}u_R^2dx\le e^{Ct}\int_{ B_R}e^{2g_0} \sigma_0^2dx,\label{eq:0order_est}\\
        &\int_{B_R}e^{2g_1}\widetilde{\nabla}u_R\cdot\nabla u_Rdx\le e^{Ct}\int_{ B_R}e^{2g_1}\widetilde{\nabla}\sigma_0\cdot\nabla \sigma_0dx+e^{Ct}Ct\int_{ B_R}e^{2g_0}\sigma_0^2dx.\label{eq:1order_est}
    \end{align}
    \end{subequations}
\end{proposition}

\medskip
For discretization purposes, we assume that the pathwise observation is given in the time sequence $0=t_0<t_1<...<t_{N_T}=T$ with equidistant $\Delta t$. Since the observation data $y_{t_j}$ is unknown before time $t_j$, the algorithm freezes $y_t$ to be $y_{t_{j-1}}$ in each time interval $[{t_{j-1}},{t_j})$. This leads to a PR-DMZ equation with observation frozen at $t_{j-1}$:
\begin{align}\label{eq:frzn_robustDMZ}
\begin{cases}
 \frac{\partial u_j}{\partial t}(x,t)+\frac{\partial}{\partial t}(h^TS^{-1})y_{t_{j-1}}u_j(x,t) \\
 \qquad\qquad= \exp{(-h^TS^{-1}y_{t_{j-1}})}(\mathcal{L}-\frac{1}{2}h^TS^{-1}h)\cdot\big(\exp{(h^TS^{-1}y_{t_{j-1}})}u_j(x,t)\big),\\
u_j(x,t_{j-1})=u_{j-1}(x,t_{j-1}),
\end{cases}
\end{align}
where $1\le j\le N_T$, $u_0(x,t_0)=\sigma_0(x)$. It is easy to check that under the invertible exponential transformation
\begin{align}\label{eq:exp_transform}
    \tilde u_j(x,t)=\exp(h^T(x,t)S^{-1}(t)y_{t_{j-1}})u_j(x,t),%\notag
\end{align}
$\tilde u_j(x,t)$ satisfies the forward Kolmogorov equation
\begin{align}\label{eq:FKEunbd}
    \frac{\partial \tilde u_j(x,t)}{\partial t}=(\mathcal{L}-\frac{1}{2}h^TS^{-1}h)\tilde u_j(x,t),
\end{align}
if and only if $u_j(x,t)$ satisfies \cref{eq:frzn_robustDMZ}. Notice that the local FKE \cref{eq:FKEunbd} is independent of the observation sequence $\{y_{t_j}\}_{j=1}^{N_T}$. Therefore, one can precompute its solver during the offline stage, that is, before any observation is given. This approach of shifting computation unrelated to observation offline is also used in other filtering methods. For example, the $S^3$ algorithm \cite{lototsky1997nonlinear} applies the Wiener chaos expansion, with the coefficients solved by a recursive system of Kolmogorov-like equations. The Yau-Yau's algorithm improves on this by solving only a single Kolmogorov equation \cref{eq:FKEunbd}.

During online stage at time $t_j,\,j\ge1$, observation data $y(t_j)$ is available to update the solution at $t_j$ as the initial condition for FKE at $[t_j,t_{j+1})$,
\begin{align}\label{eq:exptrans_frzn}
    \tilde u_{j+1}(x,t)=\exp(h^TS^{-1}(y_{t_j}-y_{t_{j-1}}))\tilde u_j(x,t).
\end{align}
The approximation solution $\hat u(x,t)$ of the PR-DMZ solution $u(x,t)$ is constructed by
\begin{align*}
    \hat u(x,t) = \sum_{j=1}^{N_T} \mathcal{X}_{[t_{j-1},\,t_j]}u_j(x,t),
\end{align*}
where $u_j(x,t)$ is derived from $\tilde u_i(x,t)$ from exponential transformation \cref{eq:exp_transform}. The approximation converges to the exact solution $u(x,t)$ in the $L^1$ sense \cite{yau2000real,luo2013complete}.

\subsection{Weighted Sobolev space}
Here we introduce the weighted Sobolev space and some of its related properties. For a function $u:\Omega\rightarrow \mathbb{R}$, where $\Omega$ is a bounded open set in $\mathbb{R}$, the derivatives of $u$ are denoted with multiindex $\alpha=(\alpha_1,...,\alpha_d)$:
\begin{align*}
    D^{\alpha}u:=\frac{\partial^{|\alpha|}u}{\partial x_1^{\alpha_1}\cdots\partial x_d^{\alpha_d}}=\partial_{x_1}^{\alpha_1}\cdots\partial_{x_d}^{\alpha_d}u,
\end{align*}
where $|\alpha|=\alpha_1+\cdots+\alpha_d,\alpha_i\ge0$.

\begin{definition}[weighted Sobolev space \cite{kufner1980weighted}]
    Let $\beta(x)$ be a nonnegative function defined on $\Omega$ and positive almost everywhere. The space $W^{k,p}(\Omega;\beta)$ is defined as the set of all functions $u=u(x)$ which are defined a.e. on $\Omega$ and satisfy
    $$\int_\Omega |D^{\alpha} u|^p\beta(x)dx< \infty,$$
    where $|\alpha|\le k$, $k\in \mathbb{N}_0$, $1\le p<\infty$. $\beta(x)$ is called the weight function.
\end{definition}  
This is a Banach space equipped with the norm 

$$||u||_{W^{k,p}(\Omega;\beta)}=\Big(\sum_{|\alpha|\le k}\int_\Omega |D^\alpha u|^p\beta(x)dx\Big)^{1/p}.$$
For $k=0$, we write $W^{0,p}(\Omega,\beta)=L^p(\Omega,\beta)$, and for $p=2$, we denote $W^{k,2}(\Omega,\beta)=H^k(\Omega,\beta)$. If $\beta(x)\equiv 1$, then the weighted space is essentially the classical Sobolev space $W^{k,p}(\Omega)$. 

As a powerful generalization of the classical Sobolev space, weighted Sobolev space has been studied for a long time. The weighted analogoues to the Sobolev embedding inequalities on bounded domains are introduced in \cite{adams2003sobolev} by Adams, and the compact embeddings on unbounded domains can be found in \cite{adams1971compact}, both with weight functions determined by $|x|$. Another type of well-known weight function is the Muckenhoupt weights. The embedding operator of Muckenhoupt $A_p$-weighted space is obtained in \cite{gol2009weighted}, for smooth domains and domains with boundary singularities. More literature on embedding operators for weighted space can be found in \cite{turesson2000nonlinear, gurka1988continuousI,opic1989continuousII,gurka1991continuousIII}.

In this paper, we consider the weight functions that vanish on the boundary of physical domain. Such weights eliminate the boundary term that arises from integration by parts, thereby simplifying the derivation of \textit{a priori} estimates in \cref{thm:prioriesti}. A natural direction is to define weights based on the distance from a point inside the domain to the boundary. Specifically, we consider the \textit{power-type weight function} introduced in \cite{kufner1980weighted}, 
\begin{align}\label{eq:power_type_weight}
    \beta(x)=(\text{dist}(x,\partial\Omega)) ^\eta,
\end{align}
where $\eta\in\mathbb{R}$, and $\text{dist}(x,\partial\Omega)$ denotes the distance from $x$ to the boundary $\partial\Omega$. For simplicity, we set $d_{\partial\Omega}(x):=\text{dist}(x,\partial\Omega)$. The associated weighted Sobolev space is denoted by $W^{k,p}(\Omega; d_{\partial\Omega},\eta)$. 

The inclusion relations of these weighted Sobolev spaces, depending on the value of $\eta$, follow directly from their definition \cite{kufner1980weighted}.
\begin{proposition}\label{prop:powerweight_embedding}
     Let $\Omega $ be a bounded domain in $\mathbb{R}^d$, then for $\eta\ge0,$ 
    \begin{align}
        W^{k,p}(\Omega)\hookrightarrow W^{k,p}(\Omega; d_{\partial\Omega}, \eta),
    \end{align}
    while for $\eta\le 0,$ 
    \begin{align}
        W^{k,p}(\Omega; d_{\partial\Omega},\eta)\hookrightarrow W^{k,p}(\Omega).
    \end{align}
\end{proposition}
This proposition entails that the weighted Sobolev space is monotonic about the exponent $\eta$.

\subsection{Tensor Train format}
Discrete tensor train format \cite{oseledets2011tensor} was proposed to compress the full array into a low-rank structure, and the decomposition can be generalized to function tensor train approximation, i.e., FTT \cite{bigoni2016spectral}. Moreover, the TT-rank obtained in DTT is used to actually represent the rank of the FTT approximation.
\paragraph{Discrete Tensor Train decomposition}
A tensor ${\bf A}\in \mathbb{R}^{n_1\times\cdots\times n_d}$ is called in TT-format if it is defined by
\begin{align*}
    {\bf A}(i_1,..,i_d) = \sum_{\widetilde{\alpha}_1=1}^{r_1}\cdots\sum_{\widetilde\alpha_{d-1}=1}^{r_{d-1}}G_1(\widetilde\alpha_0,i_1,\widetilde\alpha_1)\cdots G_d(\widetilde\alpha_{d-1},i_d,\widetilde\alpha_d),
    %{\bf A}(i_1,..,i_d) = \sum_{\bm{\alpha}\le{\bf r}}G_1(\alpha_0,i_1,\alpha_1)\cdots G_d(\alpha_{d-1},i_d,\alpha_d),
\end{align*}
where $\widetilde\alpha_0=\widetilde\alpha_d=1$, and ${\bf r}=(r_1,...,r_{d-1})$ is the TT-rank. Let $r_i\le r,n_i\le n$, then the DTT representation requires a storage cost of only $\mathcal{O}(dnr^2)$ parameters. The compression is further facilitated by QTT format. In particular, if the mode size of tensor ${\bf A}$ is a power of 2, that is, $n_i=2^L$, $i=1,..,d$. Then the tensor can be reshaped to a $dL$-dimensional one with mode size 2. The TT format of the reshaped tensor is called the QTT format of ${\bf A}$, with storage cost $\mathcal{O}(d\log_2(n)r^2)$. Computational cost is also reduced under the low-rank format. As introduced in \cite{oseledets2011tensor}, in TT format, the computational complexity of addition is $\mathcal{O}(dnr^3)$, matrix-by-vector multiplication is $\mathcal{O}(dn^2r^4)$, and the Hadamard product is $\mathcal{O}(dnr^4)$, TT-rounding is $\mathcal{O}(dnr^6)$. In the QTT format, the factor $n$ can be reduced to $\log_2(n)$. 
\paragraph{Functional Tensor Train decomposition}
The discrete tensor train decomposition naturally extends to multivariate function decomposition (e.g. \cite{gorodetsky2019continuous, strossner2024approximation}). Here we consider the continuous analogue proposed in \cite{bigoni2016spectral}. Let $u\in L^2({\bf I})$, where ${\bf I}=I_1\times I_2\times\cdots\times I_d\subset\mathbb{R}^d$ is a bounded cube. Using the theory of Hilbert-Schmidt integral operators \cite{renardy2006introduction}, $u$ admits a \textit{functional tensor train decomposition}:
\begin{align}\label{eq:FTT_decomposition}
    u = \sum_{\widetilde\alpha_1,...,\widetilde\alpha_{d-1}=1}^\infty \varphi_1(\widetilde\alpha_0;x_1;\widetilde\alpha_1) \varphi_2(\widetilde\alpha_1;x_2;\widetilde\alpha_2)\cdots\varphi_d(\widetilde\alpha_{d-1};x_d;\widetilde\alpha_d),
\end{align}
where $\widetilde\alpha_0,\widetilde\alpha_d=1$. The detailed derivation of \cref{eq:FTT_decomposition} is provided in \cref{app:FTT decomposition}. The FTT approximation results from truncating the FTT decomposition \cref{eq:FTT_decomposition}.
\begin{definition}
     [Functional Tensor Train Approximation] Let $u\in L^2(\bf{I})$. For ${\bf r}=(1,r_1,...,r_{d-1},1)$, a functional TT-rank-$\bf r$ approximation of $u$ is : 
     \begin{align}\label{eq:FTT_APPROX}
         u_{TT}({\bf x}):= \sum_{\widetilde\alpha_1,...,\widetilde\alpha_{d-1}=1}^{\bf r}\varphi_1(\widetilde\alpha_0;x_1;\widetilde\alpha_1) \varphi_2(\widetilde\alpha_1;x_2;\widetilde\alpha_2) \cdots\varphi_d(\widetilde\alpha_{d-1};x_d;\widetilde\alpha_d),
     \end{align}
     
     where $\varphi_i(\widetilde\alpha_{i-1},x_i,\widetilde\alpha_i) \in L^2(I_i)$ and $\langle\varphi_k(i,\cdot ,m),\varphi_k(i,\cdot,n) \rangle=\delta_{mn}$. $\varphi_i,(i=1,..,d)$ are called the cores of the approximation.
\end{definition}

The definition of DTT and FTT implies that the efficiency of tensor train format is determined by the TT-rank for prescribed accuracy. The following proposition reveals the relation between the regularity and the accuracy of FTT approximation. For simplicity, we restrict our attention to the case where the TT-ranks ${\bf r} = (r,...,r)$. 
\begin{proposition}[Convergence of the FTT approximation in \cite{bigoni2016spectral}] \label{prop:FTTconverg}Let $u\in H^k({\bf I})$ be a Hölder continuous function with $\gamma>1/2$ defined on the closed and bounded domain ${\bf I}\in\mathbb{R}^d$. Then the functional TT-rank-$\textbf{r}$ approximation is convergent, with residual 
    \begin{align*}
        ||u-u_{TT}||_{L^2}^2\le ||u||_{H^k({\bf I})}^2\frac{d-1}{k-1}\cdot\frac{1}{r^{k-1}}
    \end{align*}
    for $r\ge1$. Furthermore, 
    \begin{align*}
        \lim_{r\rightarrow\infty}||u-u_{TT}||_{L^2}^2=0\,\, \text{for }k>1. 
    \end{align*}
\end{proposition}
\begin{remark}
    A direct result of \cref{prop:FTTconverg} is how the complexity depends on the smoothness $k$ and the target precision. To approximate the function $u\in H^k({\bf I})$ by the truncated FTT \cref{eq:FTT_APPROX} to a prescribed $L^2$ precision $\varepsilon$, the required FTT-rank is $$r(\varepsilon) = \Big( \frac{(d-1)||u||_{H^k({\bf I})}^2}{(k-1)\varepsilon^2} \Big)^{\frac{1}{k-1}} \sim \mathcal{O}\Big( \big( \frac{1}{\varepsilon^2} \big)^{\frac{1}{k-1}}\Big).$$
\end{remark}

\Cref{prop:FTTconverg} indicates that higher regularity in the Sobolev space leads to lower FTT-rank $r$ for given accuracy, while in the context of the QTT method employed in \cref{sec:sparse approx} and \cref{sec:experiments}, the low-rank structure is established via the DTT format. Importantly, a low FTT rank implies a corresponding low DTT rank. In fact, in practical implementation of FTT construction, the FTT rank ${\bf r}$ is defined by the TT-rank of the tensor $\mathcal{V}$ formed by evaluating $u$ on a tensor grid ${\mathcal{X}}={x}_1\times\cdots\times{x}_d$ with weights ${\mathcal{W}}=w_1 \otimes\cdots\otimes w_d$ (see Procedure 1 in \cite{bigoni2016spectral}), such that
\begin{align}\label{weighted_sum}
    ||u||_{L^2({\bf I})}\approx \sum_{i_1=1}^{n_1}\cdots\sum_{i_d=1}^{n_d} u^2(\mathcal{X}_{\bf i}){\mathcal{W}_{\bf i}}
\end{align}
with prescribed accuracy. Approximating $\mathcal{V}=u(\mathcal{X})\mathcal{W}$ by a rank-${\bf r}$ DTT decomposition $\mathcal{V}_{TT}$, the DTT decomposition of $u(\mathcal{X})$ can subsequently be recovered as $\mathcal{V}_{TT}/\sqrt{\mathcal{W}}$. One can construct the rank-${\bf r}$ FTT approximation of $u$ with relative error $\varepsilon$, which arises in \cref{weighted_sum} and the DTT approximation. This result demonstrates the rationality and feasibility of applying DTT to the finite difference (FD) solution, significantly reducing the storage and computation costs.

\section{A priori estimate of solution to PR-DMZ equation}
\label{sec:a priori estimate}
Let $\Omega$ be any bounded Lipschitz domain in $\mathbb{R}^d$. Throughout the paper, we have the following assumptions. 
\begin{assumption}\label{as:ellipt_bnd}
    The matrix ${GQG}^T$ is uniformly elliptic and bounded in the cylinder $\mathbb{Q}_T=\Omega\times[0,T]$, i.e., there exists positive constants $\lambda_1\le\lambda_d$, such that 
%\begin{equation}\label{assump:bnded_coercive}
%    \lambda_1|\xi|^2\le \sum_{ij=1}^d(GQG^T)_{ij}\xi_i\xi_j\le\lambda_d |\xi|^2
%\end{equation}
\begin{align*}
    \lambda_1|\xi|^2\le \sum_{ij=1}^d(GQG^T)_{ij}\xi_i\xi_j\le\lambda_d |\xi|^2
\end{align*}
for any $(x,t)\in\mathbb{Q}_T$, $\xi=(\xi_1,...,\xi_d)^T\in\mathbb{R}^d$.
\end{assumption}
\begin{assumption}\label{as:laplace growth bnd}
    For all $t\in[0,T]$,
    \begin{align*}
        ||\frac{d}{dt}(GQG^T)||_\infty<\infty.
    \end{align*}
\end{assumption}

In real applications, we need to consider the PR-DMZ on $\Omega$ with a zero Dirichlet boundary condition: 
\begin{align}\label{eq:IBVP}
\begin{cases}
&\frac{\partial u_R}{\partial t}=\frac{1}{2}\widetilde{\Delta}u_R(x,t)+ F(x,t)\cdot\nabla u_R(x,t)+J(x,t) u_R(x,t)\\
&u_R(x,t)=0\quad\text{on }\partial\Omega\times[0,T]
\\&u_R(x,0)=\sigma_0(x)\quad \text{on } \Omega\times\{t=0\}.
\end{cases}
\end{align}
Therefore, the \textit{a priori} estimate is given for this initial boundary value problem.

In addition to this, the constant $C$ denotes a generic constant, unless explicitly stated otherwise.

\subsection{A priori estimate}
In this subsection, we derive the \textit{a priori} estimate of $s$-th order derivatives of the solution to the PR-DMZ equation \cref{eq:IBVP} on $[0,T]\times \Omega$. The \textit{a priori} estimate in \cref{prop:low_priori_est} of $ u_R $ and \( \nabla u_R \) in weighted Sobolev spaces ensures that they also lie in a classical Sobolev space. In fact, the \textit{a priori} estimate \cref{eq:low_priori_est} implies that $u_R$ lies in $L^{2}(\Omega;e^{g_0(x)})$, and $\nabla u_R$ lies in $L^{2}(\Omega;e^{g_1(x)})$. 
Then we have 
\begin{align*}
\int_{\Omega} |u_R|^2 dx &\le \int_{\Omega} e^{2g_0(x)}|u_R|^2dx < \infty,\\
\int_{\Omega} |\nabla u_R|^2 dx &\le \int_{\Omega} e^{2g_1(x)}|u_R|^2dx < \infty,
\end{align*}
i.e., 
$$u_R\in W^{1,2}(\Omega).$$

Notably, one of the keys to the completion of the proof in \cref{prop:low_priori_est} is the zero Dirichlet condition. It naturally cancel the boundary terms arises in integration by parts. However, such an auxiliary condition is not available for the \textit{a priori} estimates of higher-order derivatives. To address this problem, we replace $e^{g(x)}$ by a function with zero trace on $\partial\Omega$. 

\begin{theorem}[A priori estimate of the higher order derivatives]\label{thm:prioriesti}
     Consider the PR-DMZ equation on $\Omega\times[0,T]$. Let $C_y=\max_{0\le t\le T}|y_t|$. Assume for some integer $s\geq2$, there exists a set of nonnegative functions $\{\rho_0(x), \rho_1(x),...,\rho_{s}(x)\}$ with zero trace on $\partial\Omega $, such that for $\forall (x,t)\in\mathbb{Q}_T$,
     \begin{subequations}
         \begin{align}
            &\rho_0\le Ce^{g_0},\quad\rho_1\le Ce^{g_1},\label{as:low weight bnd}\\
            &\frac{d}{2\lambda_1} ||\frac{\partial }{\partial t} (GQG^T)|| _\infty+\lambda_1^{-1}(|\widetilde{\nabla}\rho_{m+1}|^2+\alpha!\rho_{m+1}|F|^2)+ J\le \frac{C}{2},\label{as:priori growth bnd}\\
            &{\rho_{m+1}^2}|D^\beta F| ^2 \le C{\rho_{m+1-|\beta|}^2},\,{\rho_{m+1}^2}|D^\beta J| ^2 \le C{\rho_{m+1-|\beta|}^2},\,\beta\le\alpha, 1\le|\beta|\le\alpha,\label{as:weighted bnd} \\
            & \rho_{m+1}^2|\nabla(D^{\alpha}J)|^2\le C\rho_{0}^2,  \label{as:weighted bnd J}\\
            &\int_{\Omega}\rho_0^2 \sigma_0 ^2\,dx< \infty,\int_{\Omega} {\rho_{|\beta|+1}^2} \widetilde{\nabla}  (D ^{\beta}\sigma_0) \cdot\nabla (D^{\beta}\sigma_0)dx<\infty,\beta\le\alpha, \label{as:weighted init bnd}
         \end{align}
     \end{subequations}
where $m=|\alpha|$, $1\le|\alpha|\le s-1$. Then we have the $s$-order estimates
\begin{comment}
\begin{equation}\label{eq:s-order_priori_est}
\begin{split}&\int_{\Omega} \frac{1}{2}{\rho_{m}^2} \widetilde{\nabla}(D^{m-1}u_R(t)) \cdot\nabla (D ^{m-1} u_R(t))dx\\
\le& e^{Ct}\Big(\sum_{1\le r\le m}P^{m}_{m-r}(t)\int_{\Omega} {\rho_{r}^2} \widetilde{\nabla}  (D ^{r-1}\sigma_0) \cdot\nabla (D^{r-1}\sigma_0)dx+P^{m}_{m}(t)\int_{\Omega} {\rho_{0}^2} \sigma_0^2dx\Big),
    \end{split}\end{equation}
for $t\in[0,T], 1\le m\le s$, where $P_{m-r}^{m}(t)$ denotes a polynomial of order $(m-r),0\le r\le m$, and $C$ depends on $G,Q,f,h,\mathbb{Q_T},\rho_m$ and $C_y$.
\end{comment}
\begin{equation}\label{eq:s-order_priori_est}
\begin{split}&\int_{\Omega} \frac{1}{2}{\rho_{m+1}^2} \widetilde{\nabla}(D^{\alpha}u_R(t)) \cdot\nabla (D ^{\alpha} u_R(t))dx\\
\le& e^{Ct}\Big(\sum_{\beta\le\alpha}P^{\alpha}_{\beta}(t)\int_{\Omega} {\rho_{|\beta|+1}^2} \widetilde{\nabla}  (D ^{\beta}\sigma_0) \cdot\nabla (D^{\beta}\sigma_0)dx+P^{\alpha}_{0}(t)\int_{\Omega} {\rho_{0}^2} \sigma_0^2dx\Big),
\end{split}\end{equation}
for $t\in[0,T], 0\le |\alpha|\le s-1$, where $P^{\alpha}_{\beta}(t)$ denotes a polynomial of order $(|\alpha|-|\beta|)$, $\beta\le\alpha$, $P^{\alpha}_{0}(t)$ denotes a polynomial of order $(|\alpha|+1)$, and $C$ depends on $G,Q,f,h,\mathbb{Q}_T,\rho_m$ and $C_y$.
The estimate implies that, ${D^{\alpha} u_R} \in L^{2}(\Omega;\rho_m^2), |\alpha|=m,1\le m \le s.$
\end{theorem}

\begin{proof} Direct application of assumption \cref{as:low weight bnd} to \cref{prop:low_priori_est} and the generality of constant $C$ give the estimates of $u_R$ and $\nabla u_R$ for $0\le t\le T$:
\begin{align}
    &\int_{\Omega}\rho^2_0u_R^2dx\le e^{Ct}\int_{ \Omega}\rho^{2}_0 \sigma_0^2dx,\label{eq:0order weight est}\\
    &\int_{\Omega}\rho^{2}_1\widetilde{\nabla}u_R\cdot\nabla u_R\,dx\le e^{Ct}\int_{ \Omega}\rho^{2}_1\widetilde{\nabla}\sigma_0\cdot\nabla \sigma_0\,dx+e^{Ct}Ct\int_{ \Omega}\rho^{2}_0\sigma_0^2\,dx \label{eq:1order weight est}.
\end{align}

Now we will show by induction on $s$ with the case $s=1$ being \cref{eq:1order weight est}. 

Assume that the theorem is valid for $|\alpha|=s-1$. Let $\rho_{s+1}$ be a compactly supported nonnegative function in $\Omega$ that satisfies assumptions \cref{as:priori growth bnd,as:weighted bnd,as:weighted init bnd}, and $\alpha$ be any multiindex with $|\alpha|=s$, then
\begin{equation}\notag
\begin{aligned}
    &\frac{d}{dt}\int_{\Omega}\frac{1}{2}{\rho_{s+1}^2} \widetilde{\nabla}(D ^{\alpha}u_R) \cdot\nabla (D^{\alpha}u_R)dx\\
    =&\int_{\Omega}\frac{1}{2}{\rho_{s+1}^2}(\nabla(D^{\alpha}u_R))^T\frac{d (GQG^T)}{d t}
    \nabla(D^{\alpha}u_R)+\int_ {\Omega}{\rho_{s+1}^2}\nabla (D^ {\alpha}u_R)_t\cdot \widetilde{\nabla}(D^{\alpha}u_R)dx\\
    =:&\text{I+II}\,.
\end{aligned}
\end{equation}
Due to the \cref{as:laplace growth bnd}, $I$ is bounded by 
\begin{equation}\notag
\begin{aligned}
    \text{I}&\le \frac{1}{2}||\frac {d }{d t}(GQG^T)|| _\infty\int_ {\Omega}{\rho_{s+1}^2} |\nabla(D^{\alpha}u_R)|^2dx\\
    &\le\frac{d}{2\lambda_1} ||\frac{d }{d t} (GQG^T)|| _\infty \int_{\Omega}{\rho_{s+1}^2}\nabla(D^{\alpha} u_R)\cdot \widetilde{\nabla}(D^{\alpha}u_R)dx.
\end{aligned}
\end{equation}
Apply integration by parts to $\text{II}$, 
%\begin{equation}\notag
    \begin{align*}
    \text{II}& = -\int _ {\Omega}{\rho_{s+1}}\Big(2\widetilde{\nabla}\rho_{s+1}\cdot \nabla(D^{\alpha}u_R)+\rho_{s+1} \widetilde{\Delta}(D^{\alpha}u_R)\Big)\notag
    \\&\qquad \qquad \qquad \cdot \Big(\frac{1}{2} \widetilde{\Delta}(D^{\alpha}u_R) +D^{\alpha}(F\cdot\nabla u_R ) + D ^ {\alpha} (Ju_R) \Big)dx\\
    =&  -\int _ {\Omega}{\rho_{s+1}}\Big\{\frac{1}{2}{\rho_{s+ 1}}\big( \widetilde{\Delta} (D^{\alpha}u_R)\big)^2+\Big[(\widetilde{\nabla}\rho_{s+1})\cdot\nabla(D^ {\alpha}u_R)\\
    &\underbrace{\qquad\qquad+{\rho_{s+1}}D^{\alpha}(F\cdot\nabla u_R)\Big]\widetilde{\Delta}(D ^{\alpha}u_R)\Big\}dx\qquad\qquad\quad }_{\text{III}}\\
    &\underbrace{-\int_ {\Omega}{\rho_{s+1}}\Big[(2\widetilde{\nabla}\rho_{s+1})\cdot \nabla(D^{\alpha}u_R)\Big]D^{\alpha}(F\cdot\nabla u_R )dx}_{\text{IV}}\notag
    \\&\underbrace{-\int_ {\Omega}{\rho_{s+1}}D^{\alpha}(J u_R) (2\widetilde{\nabla}\rho_{s+1} )\cdot\nabla(D^{\alpha}u_R)dx}_{\text{V}} \underbrace{- \int_{\Omega} {\rho_{s+1}^2}D^{\alpha} (Ju_R)\widetilde{\Delta} (D^{\alpha}u_R)dx}_{\text{VI}}\notag
    \\ =&: \text{III+IV+V+VI}\,,
\end{align*}
%\end{equation}
with the term $\rho_{s+1}^2$ vanishing on the boundary. 
For $\text{III+IV}$,
%\begin{equation}\notag
\begin{align*}
    \text{III+IV} =&-\int_{\Omega}\frac{1}{4}{\rho_{s+1}^2} (\widetilde{\Delta}(D^{\alpha} u_R))^2 dx\\
    &- \int_{\Omega}\frac{1}{4}\Big\{{\rho_{s+1}}\widetilde{\Delta}(D^{\alpha} u_R)+2\Big[(\widetilde{\nabla}\rho_{s+1})\cdot\nabla (D^{\alpha}u_R)+{\rho_{s+1}}D^{\alpha}(F\cdot\nabla u_R)\Big] \Big\}^2dx\\
    &+\int_{\Omega} \Big\{\Big[(\widetilde{\nabla}\rho_{s+1})\cdot\nabla(D ^{\alpha}u_R)+{\rho_{s+1}}D^{\alpha}(F\cdot\nabla u_R)\Big]^2\\
    &-2{\rho_{s+1}}\Big[(\widetilde{\nabla}\rho_{s+1})\cdot\nabla (D^{\alpha}u_R)\Big]D^{s}(F\cdot\nabla u_R ) \Big\}dx\notag
    \\\le& \int_{\Omega}\Big[(\widetilde{\nabla}\rho_{s+1}) \cdot\nabla(D^ {\alpha}u_R)\Big]^2 dx+\int_{\Omega} {\rho_{s+1}^2} \Big(\sum_{\beta\le\alpha}{\alpha \choose \beta} D ^{\alpha-\beta} F\cdot \nabla (D ^{\beta}u_R)\Big)^2 dx\notag
    \\\le& \int_{\Omega} \Big(|\widetilde{\nabla}\rho_{s+1}|^2+\alpha! \rho_{s+1}^2|F|^2\Big) | \nabla (D^{\alpha} u_R)| ^2 dx\notag
    \\&+\int _ {\Omega} \alpha!{\rho_{s+1}^2} \sum_{\beta<\alpha} |{\alpha \choose \beta}D ^{\alpha-\beta} F|^2 | \nabla (D ^\beta u_R)| ^2 dx.
\end{align*}
%\end{equation}
where ${\alpha \choose \beta}:=\frac{\alpha!}{\beta!(\alpha-\beta)!}$, $\alpha!:=\alpha_1!\cdots\alpha_d!$. With integration by parts, $\text{VI}$ turns out to be
%\begin{equation}\notag
    \begin{align*}
    \text{VI}=&\int_{\Omega} {\rho_{s+1}} D^{\alpha}(Ju_R)(2 \widetilde{\nabla}\rho_{s+1} )\cdot\nabla (D^{\alpha}u_R)dx+ \int_{\Omega}{\rho_{s+1}^2}\nabla D^{\alpha} (Ju_R) \cdot \widetilde{\nabla} (D^{\alpha} u_R)dx\notag
    \\=&-\text{V}+ \sum_{\beta\le\alpha} {\alpha \choose \beta} \Big(\int_{\Omega} {\rho_{s+1}^2} \widetilde{\nabla} (D ^{\alpha} u_R) \cdot\nabla(D^ { \alpha-\beta}J) D^\beta u_R  \\
    &+ \int_{\Omega} {\rho_{s+1}^2} \widetilde{\nabla} (D ^{\alpha} u_R) \cdot\nabla(D^\beta u_R) D^ {\alpha-\beta}J \Big)dx\notag
    \\\le& -\text{V}+ \sum_\beta {\alpha \choose \beta} \int _ {\Omega} {\rho_{s+1}^2} \Big[\frac{1}{2}\widetilde{\nabla} ( D ^{\alpha} u_R) \cdot \nabla(D^{\alpha}u_R) +\frac {\lambda_d}2 |\nabla( D ^ { \alpha-\beta}J)|^2 ( D ^\beta u_R)^2 \Big]dx\notag\\
    &+ \sum_{\beta< \alpha}{\alpha \choose \beta} \int_ {\Omega}\frac{\rho_{s+1}^2}{2}\Big[\widetilde{\nabla}(D^\beta u_R)\cdot \nabla (D^\beta u_R)(D ^ {\alpha-\beta}J)^2\\
    &+\widetilde{\nabla} ( D ^{\alpha} u_R)\cdot \nabla(D^{\alpha}u_R)\Big]dx
    +\int_{\Omega}{\rho_{s+1}^2}\widetilde{\nabla}(D^{\alpha}u_R)\cdot\nabla(D^{\alpha}u_R)Jdx\,. \notag
\end{align*}
%\end{equation}
Combine these results and by assumption \cref{as:priori growth bnd,as:weighted bnd}, we obtain
%\begin{equation}
\begin{align}\label{eq:diffgornw}
    &\frac{d}{dt}\Big(\int_{\Omega} \frac{1}{2}{\rho_{s+1}^2} \widetilde{\nabla}(D^{\alpha}u_R) \cdot\nabla(D^{\alpha} u_R)dx\Big)\\
    %\le&\frac{C}{2} \int_{\Omega} {\rho_{s+1}^2}\widetilde{\nabla}(D ^{\alpha}u_R) \cdot\nabla(D^{\alpha}u_R)dx+C\sum_{\beta<\alpha}\int_{\Omega} {\rho_{|\beta|+1}^2} \widetilde{\nabla}(D ^\beta u_R) \cdot\nabla(D ^\beta u_R)dx \notag\\
    %&+C\int_{ \Omega}{\rho_{s+1}^2}|\nabla(D^s J)|^2 u_R^2dx\notag\\
    \le &\frac{C}{2} \int_{\Omega} {\rho_{s+1}^2}\widetilde{\nabla}(D ^{\alpha}u_R) \cdot\nabla(D^{\alpha}u_R)dx+C\sum_{\beta<\alpha}\int_{\Omega} {\rho_{|\beta|+1}^2} \widetilde{\nabla}(D ^\beta u_R) \cdot\nabla(D ^\beta u_R)dx\notag\\
    &+C\int_{ \Omega}\rho_0^2u_R^2dx\,.\notag
\end{align}
%\end{equation}
Since we assume that \cref{eq:s-order_priori_est} holds true for $0\le m\le s$, we substitute it into \cref{eq:diffgornw} and have
%\begin{equation}\notag
\begin{align*}
    &\frac{d}{dt}(\int_{\Omega} \frac{1}{2}{\rho_{s+1}^2} \widetilde{\nabla}(D^{\alpha}u_R) \cdot\nabla(D^{\alpha} u_R))dx\notag\\
    \le&\frac{C}{2} \int_{\Omega} {\rho_{s+1}^2}\widetilde{\nabla}(D ^{\alpha}u_R) \cdot\nabla(D^{\alpha}u_R)dx+Ce^{Ct} \int_{ \Omega}{\rho_0^2} \sigma_0^2dx\notag
    \\&+C\sum_{\beta<\alpha}\Big(\sum_{\kappa\le\beta}P_{\kappa}^{\beta}(t)\int_{\Omega} {\rho_{|\kappa|+1}^2} \widetilde{\nabla}(D^{\kappa}\sigma_0) \cdot\nabla(D^{\kappa}\sigma_0)dx+P^{\beta}_ {0}(t)\int_{\Omega} {\rho_0^2}\sigma_0^2dx\Big)e^{Ct}\notag
    \\\le&\frac{C}{2} \int_{\Omega} {\rho_{s+1}^2}\widetilde{\nabla}(D ^{\alpha}u_R) \cdot\nabla(D^{\alpha}u_R)dx\\
    &+Ce^{Ct}\sum_{\kappa<\alpha}\tilde P_{\kappa}(t)\int_{\Omega} {\rho_{|\kappa|+1}^2} \widetilde{\nabla}  (D ^{\kappa}\sigma_0) \cdot\nabla (D^{\kappa}\sigma_0)dx+e^{Ct}\tilde P_{\alpha}(t)\int_{ \Omega}{\rho_0^2} \sigma_0^2dx,
\end{align*}
%\end{equation}
where 
\begin{equation}\notag
    \begin{aligned}
    &\tilde P_{\kappa}(t)=\sum_{\kappa\le\beta<\alpha}P^\beta_{\kappa}(t),\notag\\&\tilde P_{\alpha}(t)=\sum_{\beta<\alpha}P^\beta_{0}(t)+C.
\end{aligned}
\end{equation}
Therefore, for $\tilde P_{\kappa}(t)$ ($\kappa<\alpha$) is a polynomial of order $(|\alpha|-|\kappa|-1)$, i.e., $(s-|\kappa|-1)$, and $\tilde P_{\alpha}(t)$ is a $|\alpha|$-order polynomial. This implies
%\begin{equation}\notag
\begin{align*}
    &\frac{d}{dt}e^{-Ct}\Big(\int_{\Omega} \frac{1}{2}{\rho_{s+1}^2} \widetilde{\nabla}(D ^{\alpha}u_R) \cdot\nabla (D ^{\alpha} u_R)dx\Big)\notag
    \\=&e^{-Ct}\Big(-\frac{C}{2}\int_{\Omega} {\rho_{s+1}^2} \widetilde{\nabla} (D ^{\alpha}u_R) \cdot\nabla (D ^{\alpha} u_R)dx+ \frac{d}{dt}\int_{\Omega}\frac{1}{2}{\rho_{s+1}^2} \widetilde{\nabla}( D ^{\alpha}u_R)\cdot \nabla(D ^{\alpha} u_R)dx\Big)\notag
    \\\le& Ce^{Ct}\sum_{\kappa<\alpha}\tilde P_{\kappa}(t)\int_{\Omega} {\rho_{|\kappa|+1}^2} \widetilde{\nabla}  (D ^{\kappa}\sigma_0) \cdot\nabla (D^{\kappa}\sigma_0)dx+e^{Ct}\tilde P_{\alpha}(t)\int_{ \Omega}{\rho_0^2} \sigma_0^2dx
    %C\sum_{1\le r\le s}\tilde P_{s-r}(t)\int_{\Omega} {\rho_{r}^2} \widetilde{\nabla} (D ^{r-1}\sigma_0) \cdot\nabla (D^{r-1}\sigma_0)dx+\tilde P_{s}(t)\int_{ \Omega}{\rho_0^2} \sigma_0^2dx.
\end{align*}
%\end{equation}
Applying Gr{\" o}nwall's inequality, we then obtain \cref{eq:s-order_priori_est} for $|\alpha|=s$, that is the $(s+1)$-order estimate. %It follows directly that ${D^{\alpha} u_R} \in L^{2}(\Omega;\rho_m^2), 1\le |\alpha| \le s.$
$\hfill$
\end{proof}

Although the regularity of the DMZ equation has been studied (see, e.g., in \cite{kruse2012optimal}). The corresponding estimates for the pathwise robust DMZ equation remain unclear, particularly under zero-Dirichlet boundary condition. In \cref{thm:prioriesti}, we established the regularity estimates provided the observation path is bounded. 

The weighted integral estimate derived in \cref{thm:prioriesti} is one of several approaches to \textit{a priori} estimates (see \cite{Oleinik1962linear} for an overview). Another classical technique for second-order linear parabolic equations, proposed by S.N. Bernstein in \cite{bernstein1906generalisation}, provides pointwise interior \textit{a priori} estimates for classical derivatives of solutions via auxiliary functions. The estimates are generalized up to the boundary if $\partial\Omega$ is $C^2$, a requirement shared by the integral estimates in \cite{Oleinik1962linear}. In contrast, \cref{thm:prioriesti} avoids this regularity assumption. As shown later in \cref{cor:cube_pri_est}, our results apply even to cube domain\textemdash key for discretized algorithms.  
\begin{corollary}\label{cor:ball_pri_est}
    Consider the case where $\Omega=B_R$ (i.e., a ball centered at the origin point with radius $R$) and the weight functions $\rho_m(x),1\le m \le s$ are identical, with
    $$\rho_m(x)=\rho(x):=d_{\partial\Omega}^\eta\quad (\eta\ge 1).$$
    Under these conditions, the $s$th-order \textit{a priori} estimate \cref{thm:prioriesti} holds when $F$, $J$ and $D^{\alpha}F$, $D^{\beta}J\,(1\le|\alpha|\le s-1,1\le|\beta|\le s)$ are bounded in $[0,T]\times B_R$ and $C_y=\max_{0\le t\le T}|y_t|$ exists.
\end{corollary}
\begin{proof} With $\Omega=B_R$ and $\eta\ge1$, we have
\begin{align*}
    \rho(x)&=(R-|x|)^{\eta},x\in\Omega,\\
    \nabla\rho(x)&=-\frac{\eta x}{|x|}(R-|x|)^{\eta-1}.
\end{align*}
Then $|\nabla\rho|\leq \eta R^{\eta-1}$. Notice also that $F$, $J$ and $D^{\alpha}F$, $D^{\beta}J\,(1\le|\alpha|\le s-1,1\le|\beta|\le s)$, and that $C_y=\max_{0\le t\le T}|y_t|$ are bounded in cylinder $\mathbb{Q}_T$, with \cref{as:laplace growth bnd}, then the availability of assumption \cref{as:priori growth bnd,as:weighted bnd,as:weighted bnd J} in \cref{thm:prioriesti} is ensured. Thus, the \textit{a priori} estimates hold. $\hfill$
\end{proof}

\begin{corollary}\label{cor:cube_pri_est}
    Consider the case where $\Omega={\bf I}:=I_1\times\cdots I_d$, $I_i=(a_i,b_i)$, $a_i,b_i<\infty$ and the weight functions $\rho_m(x),1\le m \le s$ are identical, with
    $$\rho_m(x)=\rho(x):=d_{\partial\Omega}^\eta \quad \eta\ge1.$$
    If $F$, $J$ and $D^{\alpha}F$, $D^{\beta}J\,(1\le|\alpha|\le s-1,1\le|\beta|\le s)$ are bounded in $[0,T]\times B_R$, and that the maximum value $C_y=\max_{0\le t\le T}|y_t|$ exists, then the \textit{a priori} estimates \cref{thm:prioriesti} for $D ^mu_R,(1\le m\le s)$ hold.
\end{corollary}
\begin{proof} With $\Omega={\bf I}$ and $\eta\ge1$, we have
\begin{align*}
    \rho(x)&=\min_i\big( \min(x_i-a_i,b_i-x_i) \big),x\in{\bf I},
\end{align*}
being a piecewise linear function. Then the weak derivatives $\nabla\rho$ are piecewise constant and $|\nabla\rho|=1\, a.e.$. Note also that $F$, $J$ and $D^{\alpha}F$, $D^{\beta}J\,(1\le|\alpha|\le s-1,1\le|\beta|\le s)$, and $C_y=\max_{0\le t\le T}|y_t|$ are bounded in cylinder $\mathbb{Q}_T$, with \cref{as:laplace growth bnd}, the availability of assumption \cref{as:priori growth bnd,as:weighted bnd,as:weighted bnd J} in \cref{thm:prioriesti} is ensured. Thus, the \textit{a priori} estimates hold. $\hfill$
\end{proof}

\begin{remark}
    Due to the potential incompatibility between the initial condition and boundary condition, the solution to the DMZ equation may exhibit some singularity within an infinitesimal time interval $[0,t_{\varepsilon}]$, for instance, see, \cite{boyd1999compatibility,flyer2002convergence}. While with the dissipativity of the equation and the smoothing effect of the Laplace operator, the solution could attain higher order regularity after a short time $t_\varepsilon$, under the assumptions in \cref{thm:prioriesti}. %Therefore, for the remainder of this paper, we assume that the \textit{a priori }estimate holds for $t\in (t_{\varepsilon},T]$.
\end{remark}

\subsection{Improved regularity estimate}
The \textit{a priori} estimate in \cref{thm:prioriesti} demonstrates that $D^su_R,s\ge2$ belongs to a weighted Sobolev space, where the weight functions $\rho_m(x)(0\le m \le s)$ are zero on the boundary. This implies a potential blow-up behavior of the derivatives in the vicinity of the boundary. The lack of boundary regularity is resolved via weighted Sobolev embeddings. In this work, we focus on the \textit{power-type weight} defined in \cref{eq:power_type_weight}. This is motivated by two main factors: first, it serves as a typical and analytically convenient example with known embedding results. Second, its applicability is established to regular domains such as spheres and cubes (see \cref{cor:ball_pri_est,cor:cube_pri_est}).

Recall that $\Omega$ is defined as a Lipschitz domain. We have the following embedding theorem for the power-type weighted Sobolev space.
\begin{proposition}[Embedding theorem of the weighted Sobolev space in \cite{kufner1980weighted}]\label{prop:weightembed} 
    Let $\Omega$ be a Lipschitz domain, $1<p<\infty$, $0\le r<k$ and $\eta>(k-r)p-1$. Then
    $$W^{k,p}(\Omega;d_{\partial\Omega},\eta)\hookrightarrow W^{r,p}(\Omega;d_{\partial\Omega},\eta-(k-r)p).$$
\end{proposition}

If $\eta-(k-r)p\le 0$, then the weighted space embeds into a space no larger than the classical Sobolev space, as shown in \cref{prop:weightembed}. This transition is significant and of particular interest, since the classical Sobolev theory for $W^{k,p}(\Omega)$ provides sharper regularity estimates up to the boundary and improved continuity properties, under sufficient boundary smoothness and large $k$.

\begin{definition}[Hölder continuity and Hölder space]
     Suppose $\Omega\subset \mathbb{R}^d$ is an open set. A function $u:\Omega\rightarrow \mathbb{R}^d$ is said to be Hölder continuous with exponent $\gamma\in(0,1]$, if 
     \begin{equation}\notag
         |u(x)-u(y)|\le C|x-y|^\gamma.
     \end{equation}
    Define the seminorm $[u]_{0,\gamma}$ for $\gamma$-Hölder consinuous function :
    \begin{equation}\notag
        [u]_{0,\gamma}:=\sup_{x,y\in\Omega,x\neq y}\frac{|u(x)-u(y)|}{|x-y|^\gamma}.
    \end{equation}
    Denote $C^{k,\gamma}(\overline\Omega)$ as the set of all the functions $u\in C^k(\overline\Omega)$ with the norm  
    \begin{equation}\notag
        ||u||_{C^{k,\gamma}(\overline\Omega)}=\sum_{|\alpha|\le k}\sup_{x\in\overline\Omega}|D^\alpha u|+\sum_{|\alpha|=k}[D^\alpha u]_{0,\gamma}
    \end{equation}
    finite.
\end{definition}

\begin{proposition}[Extension and Sobolev embedding theorem in \cite{adams2003sobolev}]\label{prop:Extension_Embed}
    Let $u\in W^{k,p}(\Omega)$, $\Omega\subset\mathbb{R}^d$ be a bounded domain and $\partial\Omega\in C^{0,1}$. Then\\
    (1) For any bounded open set $\Omega'$ such that $\Omega\subset\subset\Omega'$, there exists a bounded linear operator
    $$E:W^{k,p}(\Omega)\rightarrow W^{k,p}(\mathbb{R}^d),$$
    such that $Eu=u$ a.e. on $\Omega$, $Eu$ has compact support in $\Omega'$. \\
    (2) If $k\ge1+\lfloor\frac{d}{2}\rfloor$, then 
    $$W^{k,2}(\Omega)\hookrightarrow C^{j,\gamma}(\overline\Omega),$$
    where $j=k-(1+\lfloor\frac{d}{2}\rfloor)$,
    \begin{align}
        \gamma\in\begin{cases}
        &(0,\frac{1}{2}],\,\text{for }d=2m-1
        \\&(0,1),\, \text{for }d=2m
        \end{cases} 
    \end{align}
    in which $m\in\mathbb{Z}^+$.
\end{proposition}

Combining the preceding theories, the following theorem states the transition from a weighted Sobolev space to a classical Sobolev space.

\begin{theorem}\label{thm:improved_regularity_est}
    Let $\Omega$ be a bounded Lipschitz domain, $u_R$ be the solution of the PR-DMZ equation on $\Omega\times[0,T]$. Assume that in \cref{thm:prioriesti} the $k$th-order a priori estimate with $\rho_m(x)=d_{\partial \Omega},(m=1,2,...,s) $ holds, i.e., $u_R\in H^{k}(\Omega, d_{\partial\Omega}, 2),\,k> d/2+2$. Then $u_R$ is Hölder continuous with $\gamma>\frac{1}{2}$ and $u_R\in H^{k-1}(\overline\Omega).$
\end{theorem}
\begin{proof}
Substitute $ p=2,\, \eta=2,r=k-1$ into \cref{prop:weightembed}, it follows that
$$W^{k,2}(\Omega;d_{\partial\Omega},2)\hookrightarrow W^{k-1,2}(\Omega).$$
Since $\partial\Omega$ is a Lipschitz boundary, by \cref{prop:Extension_Embed}, we have $u_R\in W^{k-1,2}(\overline\Omega).$

Notice that $k-1> d/2+1$. If $d$ is an even positive integer, then $u_R\in H^{m}(\Omega),m=d/2+1$. By \cref{prop:Extension_Embed}, $u_R\in C^{0,\alpha}(\overline \Omega),\alpha\in(0,1)$. If $d$ is an odd positive number, then $k-1> 1+(d+1)/2$. We have $u_R\in H^{m}(\Omega),m=(d+3)/2$. Applying \cref{prop:Extension_Embed}, we have $u_R\in C^{1,\alpha}(\overline\Omega)$. This implies that $u_R$ is Lipschitz continuous.

Therefore, $u_R\in H^{k}(\Omega,d_{\partial\Omega},2)$ implies $u_R\in W^{k-1,2}(\overline\Omega)$ and that $u_R$ is $\gamma$-Hölder continuous with $\gamma>\frac{1}{2}$. $\hfill$
\end{proof}

\section{Sparse approximation of PR-DMZ solution and its convergence}
\label{sec:sparse approx}

One of the efficient methods to solve the PR-DMZ equation is FDM FDM-based QTT method \cite{li2022solving}. The effectiveness of the method depends critically on the smoothness of the solution, as a low-rank TT/QTT format requires sufficient regularity to maintain accuracy. Furthermore, deriving the truncation error via Taylor expansion demands higher-order derivatives, reinforcing the need for a well-behaved solution. In this section, we work under the assumption that the \textit{a priori} estimate in \cref{thm:prioriesti} holds and derive the efficiency and convergence of the method in detail.

\subsection{Sparse approximation and its implementation}
\subsubsection{Sparse approximation in Tensor-Train decomposition}
By FTT format and applying \cref{prop:FTTconverg}, we prove that the regularity estimate in \cref{thm:improved_regularity_est} implies a sparse approximation.
\begin{theorem}
    Let $\Omega$ be a bounded Lipschitz domain, $u_R$ be the solution of the PR-DMZ equation on $\Omega\times[0,T]$. Assume that the $k$th-order a priori estimate in theorem \ref{thm:prioriesti} with $\rho_m(x)=d_{\partial\Omega},(m=1,...,s)$ holds, i.e.,
    $$u_R(x,t)\in H^{k}(\Omega, d_{\partial\Omega}, 2),\,k> d/2+2.$$ 
    Then $u_R$ has a sparse approximation in the FTT form $u_{R,TT}$ of rank $r$:
    \begin{align}
        u_{R,TT} = 
    \sum_{\alpha_0,...,\alpha_d=1}^{\textbf{r}}\varphi_1(\alpha_0,x_1,\alpha_1)\cdots \varphi_d(\alpha_{d-1},x_d,\alpha_d)
    \end{align}
    where $\textbf{r}=(1,r,..,r,1)$, with residual 
    \begin{align*}
        ||u_R-u_{R,TT}||_{L^2}\le ||u_R||_{H^{k-1}}\sqrt{\frac{d-1}{k-2}\cdot\frac{1}{r^{k-2}}},\,r\ge 1.
    \end{align*}
\end{theorem}
\begin{proof}
From {\cref{thm:improved_regularity_est}}, we know that $u_R\in H^{k-1}(\overline{\Omega})$ and $u_R$ is H{\"o}lder continuous with exponent $\gamma>\frac{1}{2}$. It follows from \cref{prop:FTTconverg} that $u_R$ has a rank-$r$ TT approximation and
\begin{align*}
        ||u_R-u_{R,TT}||_{L^2}\le ||u_R||_{H^{k-1}}\sqrt{\frac{d-1}{k-2}\cdot\frac{1}{r^{k-2}}},\,r\ge 1.
    \end{align*}
 $\hfill$
\end{proof}

\subsubsection{QTT-based FDM on PR-DMZ equation}
In order apply the FDM, assume that $\Omega$ is a cube: $\Omega={\bf I}=I_1\times ...\times I_d, I_i=(a_i,b_i), a_i,b_i<\infty $. The Yau-Yau's algorithm is essentially solving the local FKE repeatedly,
\begin{equation}\label{eq:FKEbd}
    \begin{aligned}
    \begin{cases}
    &\frac{\partial}{\partial t}\tilde u_j(x,t)=(\mathcal{L}-\frac{1}{2}h^TS^{-1}h)\tilde u_j(x,t)\quad \text{in } \Omega\times[t_{j-1},t_j],\\
    &\tilde u_j(x,t)=0\quad\text{on }\partial\Omega\times[t_{j-1},t_j],
    \\&\tilde u_j(x,t_{j-1})=\exp(h^TS^{-1}(y_{t_{j-1}}-y_{t_{j-2}}))\tilde u_{j-1}(x,t_{j-1})\quad \text{on } \Omega\times\{t=t_{j-1}\},
    \end{cases}
    \end{aligned}
\end{equation}
for $1\le j\le N_T$. Define $y_{t_{-1}}=0,\,y_{t_{0}}=0$ and $u_0(x,t_0)=\sigma_0(x)$. Partition the domain $\bf I$ into a uniform rectangular grid with spacing $\Delta x$, i.e.,
\begin{equation}\label{eq:spatial_partition}
    \begin{aligned}
    a_i=x_{i,0}<x_{i,1}<...<x_{i,N}<x_{i,N+1}=b_i,\,\Delta x=\frac{b_i-a_i}{N+1}. 
\end{aligned}
\end{equation}
 Similarly, we apply a uniform grid in the time direction with a spacing $\tau=\Delta t/N_\tau$. Denote the numerical approximation of $\tilde u_j(x_{\bf i}, t_{j-1}+n\tau)$ by $\tilde U_{\bf i}^{j,n}$, ${\bf i}\in\{0,1,...,N\}^d$. Applying the discretization for the time derivative, we have the semi-discrete scheme
\begin{equation}
    \begin{aligned}
        \tilde u^{j,n}(x) = \Big( \tau(\mathcal{L}-\frac{1}{2}h^TS^{-1}h)+I \Big) \tilde u^{j,n-1}(x)
    \end{aligned}
\end{equation}
for $n=1,...,N_\tau$, where $\tilde u^{j,n}(x)$ is the semi-discrete approximation of $\tilde u_j(x, t_{j-1}+n\tau)$. Subsequently, the assimilation of new observations in \cref{eq:exptrans_frzn} can be included in the semi-discrete scheme as
\begin{equation}
    \begin{aligned}
        \tilde u^{j+1,0}(x) = \exp\Big\{h^TS^{-1}(y_{t_{j}}-y_{t_{j-1}})\Big\}\tilde u^{j,N_\tau}(x).
    \end{aligned}
\end{equation}

\paragraph{procedure}
In order to implement a fully discrete scheme for problems of high dimension, the QTT format is applied to replace the full vectors and matrices that appear in FDM. Here, we impose the assumptions as follows.
\begin{assumption}\label{as:polyadic_f}
   The drift term $f$ is in functional polyadic decomposition of rank $N_f$, i.e., for $k=1,...,d$,
\begin{equation}\notag
     {f}_k = \sum_{l=1}^{N_f}{f}_{k}^{(l)}({x}) = \sum_{l=1}^{N_f} f_{k}^{(l,1)}({\bf x}_1^{(k,l)})\dotsi f_{k}^{(l,m_{(k,l)})}({\bf x}_{m_{(k,l)}}^{(k,l)}),\\ 
\end{equation}
where ${\bf x}_i^{(k,l)}\in{\bf I}_i^{(k,l)}\subset\mathbb{R}^{d_i^{(k,l)}-d_{i-1}^{(k,l)}}$ are a set of sorted elements of $x$ which are mutually disjoint, and $0=d_0^{(k,l)}<d_1^{(k,l)}<\cdots<d_{m_{(k,l)}}^{(k,l)}=d$, ${\bf I}_1\times\cdots \times{\bf I}_{m_{(k,l)}}={\bf I}$. Without loss of generality, we further assume that for $\forall \,k,l$, $m_{(k,l)}=m, d_i^{(k,l)}=d_i,(1\le i\le m)$.
\end{assumption}
\begin{assumption}\label{as:polyadic_h}
    The observation term $h$ is in functional polyadic decomposition of rank $d$, 
   \begin{equation}\notag
      {h}_k = \sum_{l=1}^{d}{h}_{k}^{(l)}({x}_l).
\end{equation}
\end{assumption}

The assumptions are based on two key arguments. From the theoretical side, the polyadic structure allows us to derive explicit complexity estimates with clarity. From the numerical implementation side, although any full array can be converted to TT/QTT format by the TT-SVD algorithm \cite{oseledets2011tensor, oseledets2010approximation}, the cost grows as $\mathcal{O}(n^dr^2)$, which becomes expensive in high dimensions. Even though the TT-cross algorithms \cite{savostyanov2011fast, savostyanov2014quasioptimality} alleviate some of this burden by approximating the TT/QTT tensor from function evaluations, practically, we found during experiments stated in \cref{sec:experiments}, under the polyadic structure, constructing TT/QTT tensor with Kronecker product is more accurate and efficient than the cross-type algorithms. 

Notably, \cref{as:polyadic_h} can be generalized to the same structure as \cref{as:polyadic_f}, allowing the \cref{alg:online_stage} to reuse the QTT construction method from \cref{alg:offline_stage}. To show the efficiency in the cubic sensor problem, we present a simplified version of this assumption to avoid redundancy.

To define the fully discrete method, denote the $d$-dimensional identity matrix with mode size $n$ as
\begin{align*}
    {\bf Id}^{(n)}_d(i_1,...,i_d;j_1,...,j_d)=
    \begin{cases}
        1 \quad i_1=j_1,...,i_d=j_d,\\
        0\quad \text{otherwise}.
    \end{cases}
\end{align*}
where $1\le i_j\le n, 1\le j\le d$. For $d=2$, we use the shorthand ${\bf Id}^{(n)}$. As an analogous notation, ${\bf 1}_d^{(n)}$ is the $d$-dimensional all-one tensor of mode size $n$, and ${\bf 1}^{(n)}$ is a shorthand for $d=1$. For simplicity, assume that $G={\bf Id}^{(d)},Q=q{\bf Id}^{(d)},S=s{\bf Id}^{(d)},\,q,s>0$. With this preparation, a spatial discretization scheme is summarized in \cref{alg:offline_stage}. 
\begin{algorithm}
\caption{Offline procedure of QTT approach to PR-DMZ}
\begin{algorithmic}[1]\label{alg:offline_stage}
    \State \{Discretize the space domain\}
    \State $x^j:=(x_{j,1},...,x_{j,N})$, $j=1,...,d$. Transform $x^j$ into QTT format.
    \State %$\mathcal{X}_1,..,\mathcal{X}_d=\text{meshgrid}(x^1,...,x^d)$\\
           $\tilde{\mathcal{X}}_{d_{i-1}+1},...,\tilde{\mathcal{X}}_{d_i}=\text{meshgrid}(x^{d_{i-1}+1},...,x^{d_i})$, $i=1,...,m$.
    \State \{Build the discretized Laplace operator by Kronecker product and summation in QTT format\}
    \State $\Delta_1=\frac{1}{(\Delta x)^2}\text{tridiag}(1,-2,1)$,
    \State $\Delta_d=\sum_{i=1}^d\Big[ \Big(\bigotimes_{1\le j\le (i-1)}{{\bf Id}^{(N)}} \Big) \otimes \Delta_1 \otimes \Big(\bigotimes_{(i+1)\le j\le d}{\bf Id}^{(N)} \Big) \Big].$
    \State \{Construct the discretized convection term in QTT format\}
    \For{$k=1\rightarrow d$}\\
    \State $\mathbf{ f }_k = \text{diag}\Big(\sum_{l=1}^{N_f}f_{k}^{(l,1)}(\tilde{\mathcal{X}}_1,..,\tilde{\mathcal{X}}_{d_1})\otimes\dotsi\otimes f_{k}^{(l,m)}(\tilde{\mathcal{X}}_{d_{m-1}+1},..,\tilde{\mathcal{X}}_d)\Big)$
    \EndFor\\
    ${\bf C}_1=\frac{1}{2\Delta x}\text{tridiag}(-1,0,1)$,\\
    ${\bf C}_d=\sum_{i=1}^d\Big[\Big(\bigotimes_{1\le j\le (i-1)}{\bf Id}^{(N)} \Big) \otimes {\bf C}_1 \otimes \Big(\bigotimes_{(i+1)\le j\le d}{\bf Id}^{(N)}\Big){\bf f}_i\Big]$.
     \State \{Construct the discretized zero-order term in QTT format\}
     \begin{equation}\notag
         \begin{aligned}
             h^Th(\mathcal{X})=&\sum_{k=1}^{d}\sum_{l<j}2\cdot{\bf 1}_{l-1}^{(N)}\otimes h_k^{(l)}(x^l)\otimes {\bf 1}_{j-l-1}^{(N)}\otimes h_k^{(j)}(x^j)\otimes {\bf 1}_{d-j}^{(N)}\\
             &+ \sum_{k=1}^{d}\sum_{l=1}^d {\bf 1}_{l-1}^{(N)}\otimes \big(h_k^{(l)}(x^l)\circ h_k^{(l)}(x^l)\big) \otimes {\bf 1}_{d-l}^{(N)}, \\
             {\bf Q}_d=&\frac{1}{2s}\text{diag}\big(h^Th(\mathcal{X})\big).
         \end{aligned}
     \end{equation}
     \State ${\bf A}={\Delta}_d-{\bf C}_d-{\bf Q}_d$.
     \State \{Compute the FKE operator $({\bf I}+\tau {\bf A})^{N_\tau}_{TT}$ in QTT format\}
\end{algorithmic}
\end{algorithm}

With all the discretized operators established in \cref{alg:offline_stage}, we form an explicit finite difference method for solving FKE \cref{eq:FKEbd}:
\begin{equation}\label{eq:FDscheme}
    \begin{aligned}
    {\bf \tilde U}^{j,n+1}=\left(\tau {\bf A}+{\bf Id}^{(N)}_d\right){\bf \tilde U}^{j,n},
\end{aligned}
\end{equation}
where ${\bf A}:=\frac{q}{2}\Delta_d-{\bf C}_d-\frac{1}{2}{\bf Q_d},{\bf \tilde U}^{j,n}=(\tilde U^{j,n}_{\bf i}),\,(0\le j\le N_T,\,0\le n\le N_\tau)$, and $q$ is the process noise. 
\begin{remark}\label{rmk:offline_complexity}
    Suppose that the QTT-rank in the offline stage is bounded by $r$. Denote $d_{\max}=\max\{d_1-d_0,...,d_m-d_{m-1}\}$. The complexity of constructing the QTT format of $f_{k}^{(l,i)}$ on the submeshgrid $( \tilde{\mathcal{X}}_{d_{i-1}+1} ,..., \tilde{\mathcal{X}} _{d_i})$ from a full array is $\mathcal{O}(N^{d_{\max}})$. Therefore, the complexity of computing $\textbf{F}_k$s by contruction of $mdN_f$ such submeshgrid and summation is $$\mathcal{O}\Big(mdN_fN^{d_{\max}}+d^2\log_2(N)r^3N_f\Big).$$ 
    Similarly, the complexity of computing $\textbf{Q}_d$ is $\mathcal{O}\Big(dN_h^2N+d^2\log_2(N)r^3N_h^2\Big)$.\\
    The matrix-by-matrix multiplication and addition are implemented to construct $C_d$ and $\Delta _d$, with a total operations of $$\mathcal{O}\Big(mdN_fN^{d_{\max}} +d^2\log_2(N)r^3N_f+d^2\log_2(N)r^4\Big).$$ 
    %Therefore, the complexity of constructing $\Delta_d,\textbf{C}_d$ and $\textbf{Q}_d$ is $$\mathcal{O}\Big(d^2N^{d_{\max}}(N_f+N_h^2) +d^2\log_2(N)r^3(N_f+N_h^2)+d^2\log_2(N)r^4\Big).$$
    TT-rounding is required in each matrix-by-matrix product to avoid the TT-rank from growing too fast. Consequently, the complexity of building $({\bf I}+\tau {\bf A})^{N_\tau}$ is $$\mathcal{O}\Big(d\log_2(N)r^6\log_2(N_\tau)\Big).$$
    Combining all these results, we have the complexity of the offline stage $$\mathcal{O}\Big(d^2\big((N^{d_{\max}}+\log_2(N)r^3)(N_f+N_h^2)+\log_2(N)r^4 \big)+d\log_2(N) r^6\log_2(\frac{\Delta T}{\tau})\Big).$$
\end{remark}

\paragraph{Online procedure}
Suppose that the drift term $f$ and observation $h$ are both independent of $t$, which enables us to compute the operator $\left(\tau {\bf A}+{\bf Id}^{(N)}_d\right)_{TT} ^{N_\tau}$ in offline stage. Additionally, the functional polyadic setting in \cref{as:polyadic_h} allows a fast construction of TT/QTT format for the exponential part in \cref{eq:exptrans_frzn}. Specifically, the meshgrid generated from the discretization of ${\bf I}$ can be constructed by the Kronecker product and the Hadamard product $\circ$.
\begin{align*}
    \mathcal{X}_1,..,\mathcal{X}_d&:=\text{meshgrid}(x^1,...,x^d),\\
    \mathcal{X}_k&=\Big(\bigotimes_{1\le j\le (k-1)}{{\bf 1}^{(N)}} \Big) \otimes x^k \otimes \Big(\bigotimes_{(k+1)\le j\le d}{\bf 1}^{(N)} \Big),\,(k=1,...,d). 
\end{align*}
Then, the discretized $\exp\Big\{h_kS^{-1}(y_{t_{j}}-y_{t_{j-1}})_k\Big\}$ is
\begin{align*}
        &\exp\{\frac{1}{s}h_k(y_{t_{j}}-y_{t_{j-1}})_k\}(\mathcal{X}_1 ,.., \mathcal{X}_{d}) \\
        =& \exp\{\frac{1}{s}\sum _{l=1}^{d}h_k^{(l)}(\mathcal{X}_l)(y_{t_{j}}-y_{t_{j-1}})_k\}\\
        =&  \exp\{\frac{1}{s}h_k^{(1)}(y_{t_{j}}-y_{t_{j-1}})_k\}(\mathcal{X}_1) \circ\cdots\circ \exp\{\frac{1}{s}h_k^{(d)}(y_{t_{j}}-y_{t_{j-1}})_k\}(\mathcal{X}_d) \\
        =&\big(\exp\{\frac{1}{s}h_k^{(1)}(y_{t_{j}}-y_{t_{j-1}})_k\}(x^1)\otimes{\bf 1}^{(N)}\otimes\cdots\otimes{\bf 1}^{(N)}\big) \circ\\
        &\cdots\circ \big({\bf 1}^{(N)}\otimes\cdots\otimes{\bf 1}^{(N)} \exp\{\frac{1}{s}h_k^{(d)}(y_{t_{j}}-y_{t_{j-1}})_k\}(x^d)\big)
        \\=&\exp\{\frac{1}{s}h_k^{(1)}(x^1)(y_{t_{j}}-y_{t_{j-1}})_k\}\otimes\cdots\otimes \exp\{\frac{1}{s}h_k^{(d)}(x^d)(y_{t_{j}}-y_{t_{j-1}})_k\}. 
\end{align*}

Then the online procedure is all about repeatedly solving the local FKE with explicit Euler scheme \cref{eq:FDscheme}, assimilating the new observation, and computing the statistics. The specific steps are stated in \cref{alg:online_stage}.
\begin{algorithm}
\caption{Online procedure of QTT approach to PR-DMZ}
\begin{algorithmic}[1]\label{alg:online_stage}
    \State Set up the initial data $\hat { U}^{1,0}_{\bf i}$ according to the distribution $\sigma_0$ of the initial state $x_0$.
    \For{$j=1\rightarrow N_T$}\do\\
        \State Solve the local FKE :
           $ \hat {\bf U}^{j,N_\tau} ={({\bf Id}_d^{(N)}+\tau {\bf A})^{N_\tau}_{TT}}{\hat {\bf U}}^{j,0}$.
        \State \{Convert the term $\exp\{\frac{1}{s}h^T(y_{t_{j}}-y_{t_{j-1}})\}$ into QTT format\}
        \For{$k=1\rightarrow d$}\do\\
        \State \{Build the QTT tensor\}
        \begin{align*}
            &\exp\{\frac{1}{s}h_k(y_{t_{j}}-y_{t_{j-1}})_k\}_{TT}(\mathcal{X}_1,..,\mathcal{X}_{d})\\
            =&\exp\{\frac{1}{s}h_k^{(1)}(x^1)(y_{t_{j}}-y_{t_{j-1}})_k\}_{TT}\otimes\cdots\otimes \exp\{\frac{1}{s}h_k^{(d)}(x^d)(y_{t_{j}}-y_{t_{j-1}})_k\}_{TT}. 
        \end{align*}
        \EndFor
        \State \{Assimilate the new observation data into the predicted solution $\hat {\bf U}^{j,N_\tau}$\}
        \begin{align*}
            \hat{\bf U}^{j+1,0}=\exp\{\frac{1}{s}h_1(y_{t_{j}}-y_{t_{j-1}})_1\}_{TT}\circ\cdots\circ\exp\{\frac{1}{s}h_d(y_{t_{j}}-y_{t_{j-1}})_d\}_{TT}\hat {\bf U}^{j,N_\tau}
        \end{align*}
        \State Scale the conditional density function with some constant $C$:$$\hat {\bf U}^{j+1,0}\gets \frac{1}{C}\hat {\bf U}^{j+1,0}$$
        \State Compute the statistics of prediction with the unnormalized density function $\hat {\bf U}^{j+1,0}$.
    \EndFor
\end{algorithmic}
\end{algorithm}

\begin{remark}
    Suppose that the QTT-rank in the online stage is bounded by $r$. The complexity of solving the local FKE and Hadamard product are both $\mathcal{O}(d\log_2(N)r^4)$. To control the growth of TT-rank, we apply TT-rounding after each Hadamard product, which requires operations of $\mathcal{O}(d\log_2(N)r^6)$. From \cref{rmk:offline_complexity} and the step stated in \cref{alg:online_stage}, the complexity of evaluating $\exp\{\frac{1}{s}h^T(y_{t_{j}}-y_{t_{j-1}})\}_{TT}$ is $\mathcal{O}(mdN_hN^{d_{\max}})$. Combining all these results, we have the computational complexity of online stage $\mathcal{O}(d^2N_hN^{d_{\max}}+d^2\log(N)r^4+d^2\log_2(N)r^6)$.
\end{remark}

\subsection{Convergence analysis}
The QTT-based FDM can be proved convergent. Specifically, the error between the QTT solution and the FD solution is governed by the prescribed accuracy applied in the QTT-format approximation and TT-rounding procedure. Consequently, if the underlying FD scheme is convergent, the QTT method also approximates exact solution, with the overall error controlled by time step, mesh size, and the chosen QTT accuracy threshold.
\subsubsection{Error estimate for the FD scheme}
The explicit Euler method \cref{eq:FDscheme} is consistent.
\begin{lemma}
    The truncation error $T_{\mathbf{i}}^{j,n}$ of the FD scheme \cref{eq:FDscheme} for the FKE \cref{eq:FKEbd} is $\mathcal{O}(\tau+(\Delta x)^2)$.
\end{lemma}
The regularity estimate \cref{thm:prioriesti} and \cref{cor:cube_pri_est} established in \cref{sec:a priori estimate} ensures that the solution admits a Taylor expansion, from which the truncation error can be deduced following the arguments in \cite{leveque1998finite}.

\begin{lemma}\label{lemma:FD_true_error_est}
    Let $\tilde U^{j,n}_{\bf i}$ denote the numerical solution obtained by FD scheme \cref{eq:FDscheme} and $\tilde u_j(x_{\bf i},t_{j-1}+n\tau)$ denote the corresponding exact solution of FKE \cref{eq:FKEbd} in $[t_{j-1},t_j]$. Assume that the drift term $f$, observation term $h$ and update part $\frac 1{s}\Delta y_{t_j}$ are bounded in $\Omega$ and $[0,T]$, i.e., $$|f|\le C_f,\quad|h|\le C_h,\quad\frac{1}{s} |h^T(y_{t_{j}}-y_{t_{j-1}})|\le C_K,$$ 
    where $s$ is the variance of observation noise ($E(dw_tdw_t^T)=s{\bf Id}^{(d)}dt$). Suppose $f$ is Lipschitz continuous, i.e.,
    \begin{align}
        |f_i(x)-f_i(x')|\le L_f|x-x'|,\quad\forall\,x,x'\in{\bf I},\,1\le i\le d.\notag
    \end{align}
    Let $\tilde{C}>dC_f$ be a constant, if the discretization parameters $(\tau,\Delta x)$ satisfies, 
    \begin{align}
        \Delta x<\frac{q}{C_f}, \quad \frac{\tau}{(\Delta x)^2}<\frac{1}{qd},\quad
        & \frac{\Delta x}{\tau} \le {\tilde C},\label{eq:stabcond}
    \end{align}
where $q$ is the variance of state noise (recall that $E(dv_tdv_t^T)=q{\bf Id}^{(d)}dt$), then we have the following estimate: 
\begin{enumerate}[(1)]
    \item For any discretization index pair $({\bf i },j)$, 
        \begin{equation}\label{eq:maxerror_single_interval}
        ||\tilde U^{j,N_\tau}_{\bf i}-\tilde u_j(x_{\bf i},t_j)||_\infty\le e^{\Delta t(dL_f+\frac{C_h}{2s})}||\tilde U^{j,0}_{\bf i}-\tilde u_j(x_{\bf i},t_{j-1})||_\infty+C'(\tau+(\Delta x)^2)
    \end{equation}
    where $C,C'$ do not depend on $\tau$ or $\Delta x$. Furthermore, at the final time $T$, 
    \begin{equation}\label{eq:maxerror_endtime}
        ||\tilde U^{N_T,\Delta t/\tau}_{\bf i}-\tilde u_{N_T}(x_{\bf i},T)||_\infty\le \frac{ C' e^{C_KN_T+(dL_f+\frac{C_h}{2s})T} } { e^{C_K+(dL_f+\frac{C_h}{2s})\Delta t}-1 }(\tau+(\Delta x)^2).
    \end{equation}
    \item The explicit FD scheme is uniformly stable, namely, for $\forall j$,
    \begin{equation}\label{eq:stab_onestep}
         \begin{aligned}
        ||\tilde U^{j,N_{\tau}}||_2&\le e^{C_K+\Delta t(\frac{2}{q} C_f\tilde C + \frac{C_h^2}{2s})}||\tilde U^{j-1,N_\tau}||_2.
        \end{aligned}
    \end{equation}
    Furthermore,
    \begin{equation}\label{eq:stab_general}
        \begin{aligned}
    ||\tilde U^{N_T,N_{\tau}}||_2 \le& e^{C_KN_T+T(\frac{2}{q} C_f\tilde C + \frac{C_h^2}{2s})}||\tilde U^{1,0}||_2.
\end{aligned}
    \end{equation}
\end{enumerate}
\end{lemma}
Before starting the proof, here we list two remarks.
\begin{remark} 
 The region for mesh size and time step $(\Delta x,\tau)$ is not an empty set. In particular, the stability conditions \cref{eq:stabcond} is equivalent to 
    \begin{equation}\notag%\label{eq:stab_region}
 \frac{qd}{\tilde C}<\Delta x<\frac{q}{C_f}, \quad
            \frac{\Delta x}{\tilde C}<\tau<\frac{(\Delta x)^2}{qd}.
    \end{equation}
    By $\tilde C>dC_f$, we have $ {qd}/{\tilde C}<\frac{q}{C_f}$,
    which ensures the existence of a valid mesh size $\Delta x$.
  \end{remark}
  \begin{remark}
The factor $e^{C_KN_T}$ in the estimate can be controlled under a high probability. Specifically, for any $m>0$ and $n<\frac{C_h\sqrt{\Delta t}}{m\sqrt{ds}},$ the following inequality holds true with probability $(\Phi(n)-\Phi(-n))^{d}$:
\begin{align*}
    |y_{t_{j}}-y_{t_{j-1}}|\le C_h\Delta t+n\sqrt{ds\Delta t}
        \le  (1+\frac{1}{m})C_h\Delta t.
\end{align*}
where 
\begin{align*}
    \Phi(x):=\int_{-\infty}^x \frac{1}{\sqrt{2\pi}}e^{-\frac{x^2}{2}}dx.
\end{align*}
Consequently, if the ratio $\frac{C_h\sqrt{\Delta t}}{m\sqrt{ds}}$ is sufficiently high, then $e^{\frac{1}{s}h^T(y_{t_{j}}-y_{t_{j-1}})}$ is bounded by $e^{(1+\frac{1}{m})C_h\Delta t}$ with high probability, which implies $e^{C_KN_T}\le e^{(1+\frac{1}{m})C_h T}$.
\end{remark}

\textbf{\textit{Proof of \cref{lemma:FD_true_error_est}.}}

(\textbf{Step 1}) Let $E_{\bf i}^{j,n}:=\tilde U^{j,n}_{\bf i}-\tilde u_j(x_{\bf i},t_{j-1}+n\tau)$ denote the point-wise error between the FD solution and exact solution of the FKE \cref{eq:FKEbd}. By the definition of truncation error, we have
\begin{align*}
    \frac{E_{\bf i}^{j,n+1}-E_{\bf i}^{j,n}}{\tau} 
    =&  T_{\textbf{i}}^{j,n} + \frac{q}{2(\Delta x)^2} \Big[ \sum_{k=1}^d(E^{j,n}_{{\bf i}-{\bf e}_k}+E^{j,n}_{{\bf i}+{\bf e}_k})-2dE^{j,n}_{{\bf i}} \Big] -\frac{h^Th}{2s}(x_{{\bf i}})E_{\bf i}^{j,n}\\
    - \frac{1}{2\Delta x}& \Big[\sum_{k=1}^d f_k(x_{\bf i})(E^{j,n}_{{\bf i}-{\bf e}_k} - E^{j,n}_{{\bf i}+{\bf e}_k}) + {E_{\bf i}^{j,n}}\sum_{k=1}^d \big(f_k(x_{{\bf i}+{\bf e}_k}) - f_k(x_{{\bf i}-{\bf e}_k})\big)\Big],
\end{align*}
i.e.,
\begin{align*}
    E_{\bf i}^{j,n+1} 
    =&  \tau T_{\textbf{i}}^{j,n} + \Big(1-\frac{qd\tau}{(\Delta x)^2}-\frac {\tau\sum_{k=1}^d (f_k(x_{{\bf i}+{\bf e}_k}) - f_k(x_{{\bf i}-{\bf e}_k}))} {2\Delta x} - \tau\frac{h^Th}{2s}(x_{{\bf i}}) \Big) E_{\bf i}^{j,n} \\
    &+ (\frac{qd\tau}{2(\Delta x)^2} - \sum_{k=1}^d \frac{\tau f_k(x_{\bf i})}{2\Delta x})E^{j,n}_{{\bf i}+{\bf e}_k} + (\frac{qd\tau}{2(\Delta x)^2} + \sum_{k=1}^d \frac{\tau f_k(x_{\bf i})}{2\Delta x})E^{j,n}_{{\bf i}-{\bf e}_k}.
\end{align*}
where $T_{\textbf{i}}^{j,n}$ is the truncation error. Denote the maximum error at time $t=t_j+n\tau$ as
\begin{align*}
    \mathcal{E}^{j,n}:=\max\{|E_{\bf i}^{j,n}|,{\bf i}=(i_1,...,i_d),1\le i_l\le N,l=1,...,d\}.
\end{align*}
With the stability condition \eqref{eq:stabcond}, we have
\begin{align*}
    |\frac{\tau f_k(x)}{2(\Delta x)}| \le & \frac{\tau }{2(\Delta x)}C_f \le  \frac{\tau}{2(\Delta x)}\frac{q}{(\Delta x)}=  \frac{\tau q} {2(\Delta x)^2}.
\end{align*}
Provided that $f$ is Lipschitz continuous, the inequality follows:
\begin{align*}
    (f_k(x_{{\bf i}+{\bf e}_k}) - f_k(x_{{\bf i}-{\bf e}_k}))\le 2L_f\Delta x.
\end{align*}
Therefore,
\begin{align*}
    \sum_{k=1}^d (f_k(x_{{\bf i}+{\bf e}_k}) - f_k(x_{{\bf i}-{\bf e}_k}))\le 2dL_f\Delta x.
\end{align*}
Thus, we can bound the maximum error $\mathcal{E}^{j,n+1}$ with respect to $\mathcal{E}^{j,n}$,
\begin{equation}\label{maxerror_onestep}
    \begin{aligned}
    \mathcal{E}^{j,n+1}\le& \tau T_{\textbf{i}}^{j,n} + (1-\frac{qd\tau}{(\Delta x)^2} + d\tau L_f + \frac{\tau C_h}{2s}) \mathcal{E}^{j,n} \\
    & +(\frac{qd\tau}{2(\Delta x)^2} - \sum_{k=1}^d \frac{\tau f_k(x_{\bf i})}{2\Delta x}) \mathcal{E}^{j,n} +  (\frac{qd\tau}{2(\Delta x)^2} + \sum_{k=1}^d \frac{\tau f_k(x_{\bf i})}{2\Delta x}) \mathcal{E}^{j,n} \\
    =&\tau C'(\tau+(\Delta x)^2) + \left(1 + \tau(dL_f+\frac{C_h}{2s})\right)\mathcal{E}^{j,n},
\end{aligned}
\end{equation}
where $C'$ is a generic constant. Applying \cref{maxerror_onestep} recursively to the local FKE in $[t_{j-1},t_j]$, we obtain the error estimate at time $t_j$ in terms of the error at initial time $t_{j-1}$,
\begin{equation}\label{maxerror_single_interval}
    \begin{aligned}
    \mathcal{E}^{j,N_\tau}\le e^{(dL_f+\frac{C_h}{2s})\Delta t}\mathcal{E}^{j,0} + {C'}(\tau+(\Delta x)^2),
\end{aligned}
\end{equation}
which gives \cref{eq:maxerror_single_interval}. Notice that the initial data for FKE in $[t_{j-1},t_j]$ is obtained by including the new observation $y(t_{j-1})$ into the predicted solution $\tilde u_{j-1}(x_{\bf i},t_{j-1})$. Thus,
\begin{equation}\label{error_initdata}
    \begin{aligned}
    \mathcal{E}^{j,0}\le \exp(\frac{1}{s}h^T(y_{t_{j-1}}-y_{t_{j-2}}))\mathcal{E}^{j-1,N_\tau}\le e^{C_K}\mathcal{E}^{j-1,N_\tau}.
\end{aligned}
\end{equation}
Combining \cref{maxerror_single_interval} and \cref{error_initdata}, 
\begin{align}\label{maxerror_single_interval_exp}
    \mathcal{E}^{j,N_\tau}\le e^{C_K+(dL_f+\frac{C_h}{2s})\Delta t}\mathcal{E}^{j-1,N_\tau} + C'(\tau+(\Delta x)^2).
\end{align}
Since $\mathcal{E}^{1,0}=0$, the error bound at final time follows by applying \cref{maxerror_single_interval_exp} recursively,
\begin{align*}
    \mathcal{E}^{N_T,N_\tau}\le& e^{(C_K+(dL_f+\frac{C_h}{2s})\Delta t)(N_T-1)}(\mathcal{E}^{1,N_\tau}+\frac{C'(\tau+(\Delta x)^2)}{e^{C_K+(dL_f+\frac{C_h}{2s})\Delta t}-1})\notag\\
    \le & e^{(C_K+(dL_f+\frac{C_h}{2s})\Delta t)(N_T-1)} \\
    &\cdot\left(e^{(dL_f+\frac{C_h}{2s})\Delta t}\mathcal{E}^{1,0}+C'(\tau+(\Delta x)^2)(1+\frac{1}{e^{C_K+(dL_f+\frac{C_h}{2s})\Delta t}-1})\right)\notag\\
    \le & \frac{ C' e^{C_KN_T+(dL_f+\frac{C_h}{2s})T} } { e^{C_K+(dL_f+\frac{C_h}{2s})\Delta t}-1 }(\tau+(\Delta x)^2).
\end{align*}

(\textbf{Step 2})  From the derived FD scheme, $\Delta_d$, ${\bf F}_k$ and ${{\bf Q}_d}$ are symmetric. To give the energy estimate of the FD solution, we first estimate the 2-norm of $ (\tau(\frac{q}{2}\Delta_d) +{\bf Id}^{(N)}_d)$. With the fact that $\Delta _1$ is negative definite and 
\begin{align*}
    \min(\lambda(\Delta _1)) \ge  -\frac{4}{(\Delta x)^2},
\end{align*}
combined this with the stability condition \cref{eq:stabcond}, we have the minimum eigenvalue estimate 
\begin{align*}
    \min\Big(\lambda\big(\tau(\frac{q}{2}\Delta_d) +{\bf Id}^{(N)}_d \big) \Big)&\ge 1+\frac{dq\tau}{2}\min(\lambda(\Delta _1))\\
    &\ge1 -\frac{2dq\tau}{(\Delta x)^2} \ge -1.
\end{align*}
Therefore, $\lambda( \tau(\frac{q}{2}\Delta_d) +{\bf Id}^{(N)}_d) \in [-1,1]$. Further, by the Gershgorin circle theorem, we have
\begin{align*}
    || \tau(\frac{q}{2}\Delta_d-{\bf C}_d-\frac{1}{2}{\bf Q_d})+{\bf Id}^{(N)}_d||_2 
   \le & || \tau(\frac{q}{2}\Delta_d) +{\bf Id}^{(N)}_d ||_2 + ||\tau ({\bf C}_d-\frac{1}{2}{\bf Q}_d)||_2 \notag\\
    \le& 1+ \tau||{\bf C}_d||_2+\frac{\tau}{2}||{\bf Q}_d||_2  \notag\\
    \le& 1 +  \tau (\frac{2C_fd}{\Delta x}+\frac{C_h^2}{2s})\notag\\
    \le & 1 + \tau ( \frac{2C_f\tilde C}{q} + \frac{C_h^2}{2s}).
\end{align*}
The last inequality follows from \cref{eq:stabcond},
\begin{align*}
    \frac{1}{\Delta x}=&\frac{\tau}{(\Delta x)^2}\cdot\frac{\Delta x}{\tau}
    \le \frac{\tilde C}{qd}.
\end{align*}
Therefore, 
\begin{align}\label{stab_onestep}
    ||\tilde {\bf U}^{j,N_{\tau}}||_2 &= || (\tau{\bf A}+{\bf Id}_d^{(N)})^{N_\tau}  \tilde {\bf U}^{j,0}||_2 \\
    &=||( \tau(\frac{q}{2}\Delta_d-{\bf C}_d-\frac{1}{2}{\bf Q_d})+{\bf Id}^{(N)}_d)^{N_\tau} \tilde  {\bf U}^{j,0}||_2\notag\\
    &\le (1 + \tau (  \frac{2}{q} C_f\tilde C + \frac{C_h^2}{2s}))^{N_\tau}||\tilde  {\bf U}^{j,0}||_2\notag\\
    &\le e^{\Delta t( \frac{2C_f\tilde C}{q}  + \frac{C_h^2}{2s})}||\exp\{\frac{1}{s}h^T(y_{t_{j-1}}-y_{t_{j-2}})\}\circ \tilde  {\bf U}^{j-1,N_\tau}||_2\notag\\
    &\le e^{C_K + \Delta t( \frac{2}{q}C_f\tilde C  + \frac{C_h^2}{2s})}||\tilde  {\bf U}^{j-1,N_\tau}||_2.\notag
\end{align}
Applying \eqref{stab_onestep} recursively, it leads to
\begin{align*}
    ||\tilde  {\bf U}^{N_T,N_{\tau}}||_2 \le& e^{(C_K + \Delta t(\frac{2}{q} C_f\tilde C + \frac{C_h^2}{2s}))(N_T-1)}||\tilde  {\bf U}^{1,N_\tau}||_2\notag\\
    \le& e^{(C_K + \Delta t(\frac{2}{q} C_f\tilde C + \frac{C_h^2}{2s}))N_T}||\tilde  {\bf U}^{1,0}||_2\notag\\
    = & e^{C_KN_T+T(\frac{2}{q} C_f\tilde C + \frac{C_h^2}{2s})}||\tilde  {\bf U}^{1,0}||_2.
\end{align*}
$\hfill\square$

\subsubsection{Error estimate for QTT method}
Assume that the given accuracy $\varepsilon_1$ is used in constructing the QTT-format and TT-rounding, and $\varepsilon_2$ is the error between the QTT-format discretized operator and the FD discretized operator, i.e., 
\begin{align*}
    || \hat{ \bf U }^{1,0} - \tilde {\bf U}^{1,0} ||_2 &\le \varepsilon_1 || \tilde {\bf U}^{1,0}||_2,  \\
    \Big|\Big| { \exp \{\frac{1}{s}h^T(y_{t_{j}}-y_{t_{j-1}})\}_{TT}  }  
    - \exp \{\frac{1}{s}h^T(y_{t_{j}}-y_{t_{j-1}})\}  \Big|\Big|_2
    & \le \varepsilon_1  C_{\exp},\\
    \big|\big| ( { \bf Id  } + \tau { \bf A } ) ^ { N_\tau } _{TT}
   - ( { \bf Id  } + \tau { \bf A } ) ^ { N_\tau } \big|\big|_2 
   &\le \varepsilon_2 C_{\bf A},
\end{align*}
where $C_{\exp} = \max_j\big|\big|  \exp \{\frac{1}{s}h^T(y_{t_{j}}-y_{t_{j-1}})\}   \big|\big|_2$, $C_{\bf A} = || ( { \bf Id  } + \tau { \bf A } ) ^ { N_\tau } ||_2 $.

\begin{lemma}\label{lemma:QTT_FD_error_est}
    Let $\hat {U}^{j,n}_{\bf i}$ be the QTT solution of $\tilde u(x_{\bf i}, t_j+n\tau)$, $j=1,...,N_T,n=0,...,N_\tau$, where
    \begin{align*}
        \hat {\bf U}^{j,0} = { \exp \{\frac{1}{s}h^T(y_{t_{j}}-y_{t_{j-1}})\}_{TT}  } \circ \hat {\bf U}^{j-1,N_\tau}.
    \end{align*}
    Denote $\varepsilon_0 = \max\{\varepsilon_1,\varepsilon_2\}$, then
    \begin{align}\label{error_final_QTT}
        || \hat{ \bf U }^{N_T, N_\tau} - \tilde
        { \bf U }^{N_T, N_\tau} ||_2 \le \varepsilon_0 C ||\tilde {\bf U}^ {1,0}||_2,
    \end{align}
    where
    \begin{align}
    & C=  C_1^{N_T-1}( e^{\Delta t(\frac{2}{q} C_f\tilde C + \frac{C_h^2}{2s})} + \varepsilon_2 C_{\bf A}) +   C_{\bf A}C_1^{N_T-1} + C_2\frac{C_1^{N_T} - e^{C_KN_T+T(\frac{2}{q} C_f\tilde C + \frac{C_h^2}{2s})}}{C_1e^{-\Delta t(\frac{2}{q} C_f\tilde C + \frac{C_h^2}{2s})}-1}, \label{constant_QTT_FD_error}
    \end{align}
    with $C_1 = (e^{\Delta t(\frac{2}{q} C_f\tilde C + \frac{C_h^2}{2s})} + \varepsilon_2C_{\bf A} )
  (e^{C_K} + \varepsilon_1 C_{\exp})$, $C_2 = (e^{\Delta t(\frac{2}{q} C_f\tilde C+ \frac{C_h^2}{2s})} + \varepsilon_2C_{\bf A} )C_{\exp} +C_{\bf A}e^{C_K}$.
\end{lemma}
\begin{proof}
First, consider the error generated in the step of assimilating new observation data in QTT format.
\begin{align}\label{exp_error}
    &|| \hat{ \bf U }^{j+1,0} - \tilde {\bf U }^{j+1,0} ||_2 \notag\\
    = & \Big|\Big| \exp \{\frac{1}{s}h^T(y_{t_{j}}-y_{t_{j-1}})\}_{TT}\circ \hat {\bf U}^{j,N_\tau} - \exp \{\frac{1}{s}h^T(y_{t_{j}}-y_{t_{j-1}})\}\circ \tilde {\bf U}^{j,N_\tau} \Big|\Big|_2 \notag\\
    \le & \big|\big| \exp \{\frac{1}{s}h^T(y_{t_{j}}-y_{t_{j-1}})\} \big|\big|_\infty ||\hat {\bf U}^{j,N_\tau} - \tilde {\bf U}^{j,N_\tau}||_2 \notag\\
    & + \big|\big| \exp \{\frac{1}{s}h^T(y_{t_{j}}-y_{t_{j-1}})\}_{TT} - \exp \{\frac{1}{s}h^T(y_{t_{j}}-y_{t_{j-1}})\} \big|\big|_\infty || \hat {\bf U}^{j,N_\tau} ||_2\notag\\
    \le & e^{C_K} ||\hat {\bf U}^{j,N_\tau} - \tilde {\bf U}^{j,N_\tau}||_2 + \varepsilon_1 C_{\exp} (||\hat {\bf U}^{j,N_\tau} - \tilde {\bf U}^{j,N_\tau}||_2 +  ||\tilde {\bf U}^{j,N_\tau}||_2)\notag\\
    \le & (e^{C_K} + \varepsilon_1 C_{\exp})||\hat {\bf U}^{j,N_\tau} - \tilde {\bf U}^{j,N_\tau}||_2 + \varepsilon_1 C_{\exp} e^{(C_K+\Delta t(\frac{2}{q} C_f\tilde C + \frac{C_h^2}{2s}))j}||\tilde {\bf U}^{1,0}||_2.
\end{align}
With \cref{exp_error} and the energy estimate \cref{eq:maxerror_single_interval}, we can give the error estimate between the QTT solution $\hat {\bf U}^{j+1,N_\tau}$ and FD solution $\tilde {\bf U}^{j+1,N_\tau}$, with respect to the error generated in last time interval $[t_{j-1},t_j)$.
%\begin{equation}
    \begin{align}\label{error_FKE_local}
    &|| \tilde {\bf U }^{j+1,N_\tau} - \hat {\bf U }^{j+1,N_\tau} ||_2 \\
= & \Big|\Big| ({\bf Id} + \tau {\bf A})^{N_\tau}_{TT} \hat { \bf U }^{j+1,0}
    - ({\bf Id} + \tau {\bf A})^{N_\tau} \tilde { \bf U }^{j+1,0} \Big|\Big|_2 \notag\\
\le & || ({\bf Id} + \tau {\bf A})^{N_\tau}||_2 || \hat{ \bf U }^{j+1,0} - \tilde {\bf U }^{j+1,0} ||_2 + \big|\big| ({\bf Id} + \tau {\bf A})^{N_\tau}_{TT} -  ({\bf Id} + \tau {\bf A})^{N_\tau} \big|\big|_2 || \hat{ \bf U }^{j+1,0}||_2 \notag\\
\le & e^{\Delta t(\frac{2}{q} C_f\tilde C + \frac{C_h^2}{2s})}|| \hat{ \bf U }^{j+1,0} - \tilde {\bf U }^{j+1,0} ||_2 + \varepsilon_2 C_{\bf A}(|| \hat{ \bf U }^{j+1,0} - \tilde {\bf U }^{j+1,0} ||_2 + ||\tilde {\bf U}^ {j+1,0}||_2) \notag\\
\le & (e^{\Delta t(\frac{2}{q} C_f\tilde C + \frac{C_h^2}{2s})} + \varepsilon_2C_{\bf A} )\cdot\notag\\
&\left(  (e^{C_K} + \varepsilon_1 C_{\exp})||\hat {\bf U}^{j,N_\tau} - \tilde {\bf U}^{j,N_\tau}||_2 + \varepsilon_1 C_{\exp} e^{(C_K+\Delta t(\frac{2}{q} C_f\tilde C + \frac{C_h^2}{2s}))j} ||\tilde {\bf U}^{1,0}||_2 \right)  \notag\\
& + \varepsilon_2C_{\bf A}e^{C_K} e^{(C_K+\Delta t(\frac{2}{q} C_f\tilde C + \frac{C_h^2}{2s}))j}||\tilde {\bf U}^{1,0}||_2 \notag\\
\le & C_1 ||\hat {\bf U}^{j,N_\tau} - \tilde {\bf U}^{j,N_\tau}||_2  
+  \varepsilon_0  C_2e^{(C_K+\Delta t(\frac{2}{q} C_f\tilde C + \frac{C_h^2}{2s}))j} ||\tilde {\bf U}^{1,0}||_2.\notag
\end{align}
%\end{equation}
By applying \eqref{error_FKE_local} recursively, we have the following error estiamte
%\begin{equation}\notag
    \begin{align*}
   & || \tilde {\bf U }^{N_T,N_\tau} - \hat {\bf U }^{N_T,N_\tau} ||_2 \\
   \le & C_1^{N_T-1} || \tilde {\bf U }^{1,N_\tau} - \hat {\bf U }^{1,N_\tau} ||_2 \\
   &+ \varepsilon_0 C_2 e^{C_KN_T+T(\frac{2}{q} C_f\tilde C + \frac{C_h^2}{2s})}||\tilde {\bf U}^{1,0}||_2\sum_{k=0}^{N_T-1}C_1^ke^{-(C_K+\Delta t(\frac{2}{q} C_f\tilde C+ \frac{C_h^2}{2s}))k} \\
   \le & C_1^{N_T-1} \left( e^{\Delta t(\frac{2}{q} C_f\tilde C+ \frac{C_h^2}{2s})}|| \hat{ \bf U }^{1,0} - \tilde {\bf U }^{1,0} ||_2 + \varepsilon_2 C_{\bf A}(|| \hat{ \bf U }^{1,0} - \tilde {\bf U }^{1,0} ||_2 + ||\tilde {\bf U}^ {1,0}||_2) \right) \\
   & + \varepsilon_0 C_2 e^{C_KN_T+T(\frac{2}{q} C_f\tilde C + \frac{C_h^2}{2s})}||\tilde {\bf U}^{1,0}||_2 
   \frac{1-C_1^{N_T}e^{-(C_K+\Delta t(\frac{2}{q} C_f\tilde C + \frac{C_h^2}{2s})){N_T}}}{1-C_1e^{-(C_K+\frac{2\Delta t}{q} C_f\tilde C+\frac{C_h^2\Delta t}{2s})}}\\
   \le & \varepsilon_1 C_1^{N_T-1}( e^{\Delta t(\frac{2}{q} C_f\tilde C + \frac{C_h^2}{2s})} + \varepsilon_2 C_{\bf A})||\tilde {\bf U}^{1,0}||_2 \\
   & +   \varepsilon_0 \left( C_{\bf A}C_1^{N_T-1} + C_2\frac{C_1^{N_T} - e^{C_KN_T+T(\frac{2}{q} C_f\tilde C + \frac{C_h^2}{2s})}}{C_1e^{\Delta t(\frac{2}{q} C_f\tilde C+ \frac{C_h^2}{2s})}-1}\right)||\tilde {\bf U}^{1,0}||_2 \notag\\
   \le & \varepsilon_0  ||\tilde {\bf U}^ {1,0}||_2\Big( C_1^{N_T-1}( e^{\Delta t(\frac{2}{q} C_f\tilde C + \frac{C_h^2}{2s})} + \varepsilon_2 C_{\bf A})  \notag\\
   & + C_{\bf A}C_1^{N_T-1} + C_2\frac{C_1^{N_T} - e^{C_KN_T+T(\frac{2}{q} C_f\tilde C+ \frac{C_h^2}{2s})}}{C_1e^{-\Delta t(\frac{2}{q} C_f\tilde C + \frac{C_h^2}{2s})}-1}\Big) .
\end{align*}
%\end{equation}
$\hfill$
\end{proof}

As a direct combination of \cref{lemma:FD_true_error_est} and \cref{lemma:QTT_FD_error_est}, the error estimate between the QTT solution $\hat U^{j,n}_{\bf i}$ and the exact solution $\tilde u_j(x_{\bf i},t_{j-1}+n\tau)$ is immediately available, as stated in the following theorem.
\begin{theorem}\label{thm:QTT_true_error_est}
     Let $\hat {U}^{j,N_\tau}_{\bf i}$ be the QTT solution of $\tilde u_j(x_{\bf i}, t_{j-1}+n\tau)$, $j=1,...,N_T$. Suppose the assumptions in \cref{lemma:FD_true_error_est} and \cref{lemma:QTT_FD_error_est} are all satisfied, then the maximum error estimate of $\hat {U}^{j,N_\tau}_{\bf i}$ is 
     \begin{align}
         ||\hat{U}^{N_T,N_\tau}_{\bf i}-\tilde u_j(x_{\bf i},T)||_{\infty} \le \overline{C}\left(\varepsilon_0+\tau+(\Delta x)^2\right),
     \end{align}
     where 
     \begin{align*}
         \overline C=\max\{C||\tilde {\bf U}^{1,0}||_2,\frac{ C' e^{C_KN_T+(dL_f+\frac{C_h}{2s})T} } { e^{C_K+(dL_f+\frac{C_h}{2s})\Delta t}-1 }\},
     \end{align*}
     $C$ is defined in \cref{constant_QTT_FD_error}, $C'$ is the generic constant in \cref{maxerror_onestep}.
\end{theorem}
\begin{remark}
Theorem \cref{lemma:QTT_FD_error_est} implies that, by choosing proper $\varepsilon_1,\varepsilon_2$, the error of the QTT solution can be controlled to match the error scale of FD scheme.
\end{remark}
\begin{proof}
Simply apply triangle inequality to \cref{lemma:FD_true_error_est} and \cref{lemma:QTT_FD_error_est}, we have
\begin{align}
    ||\hat{U}^{N_T,N_\tau}_{\bf i}-\tilde u_j(x_{\bf i},T)||_{\infty} &\le ||\hat{U}^{N_T,N_\tau}_{\bf i} - \tilde{U}^{N_T,N_\tau}_{\bf i}||_\infty + ||\tilde{U}^{N_T,N_\tau}_{\bf i}-\tilde u_j(x_{\bf i},T)||_{\infty}\notag \\
    &\le \varepsilon_0 C ||\tilde {\bf U}^ {1,0}||_2 + \frac{ C' e^{C_KN_T+(dL_f+\frac{C_h}{2s})T} } { e^{C_K+(dL_f+\frac{C_h}{2s})\Delta t}-1 }(\tau+(\Delta x)^2)\notag\\
    &\le \overline{C}\left(\varepsilon_0+\tau+(\Delta x)^2\right).
\end{align}
$\hfill$
\end{proof}

\section{Numerical Experiments}\label{sec:experiments}
In this section, we present two synthetic numerical examples to illustrate the performance of the proposed QTT approach to PR-DMZ (offline  \cref{alg:offline_stage} and online \cref{alg:online_stage}). The first example shows the performance of QTT in 10-20 dimensional cubic sensor problems and comparison with linear filters (EKF) proposed in \cite{jazwinski2013stochastic} and particle filters (PF) in \cite{gordon1993novel}. In the second example, we consider a multi-mode problem where EKF and PF have difficulties to track the solutions accurately. The relative mean square error (RMSE) is applied to compare the performance of different algorithms:
\begin{align*}
    \text{RMSE}:=\frac{1}{dN_T}\sum_{j=1}^{N_T} \big| \hat x(t_j) -  x(t_j) \big|^2,
\end{align*}
where $x(t_j)$ is the real state at time $t_j$ and $\hat x(t_j)$ is the estimation of $x(t_j)$. The numerical experiments were performed in MATLAB R2024b on a single node of HKU's HPC cluster, utilizing an 8-core CPU.
\subsection{Cubic sensor}
The cubic sensor \cite{hazewinkel1983nonexistence} problem considers the following nonlinear observation model:
\begin{equation}\label{eq:cubicsensor}
    \begin{aligned}
\begin{cases}
&d x_t=(B_dx_t+\sin (x_t))dt+dv_t\\
&dy_t= h(x_t)dt+dw_t\\
\end{cases}
\end{aligned}
\end{equation}
where $E(dv_tdv_t^T) = 0.8I_ddt, E(dw_tdw_t^T)=2 I_d$, and $B_n=(b_{ij})_{d\times d}$ is a square matrix with elements
\begin{align}
b_{ij}=
\begin{cases}
-0.6, &i=j,\\
-0.1,  &j-i=1,\\
0, &\text{otherwise.}\\
\end{cases}\label{Bn}
\end{align}
The observation function is $h(x)=[x_1^3,...,x_d^3]^T$. The initial unnormalized density function is $\sigma_0(x)=e^{-5{|x|^4}}$. The time step to generate a real trajectory and observation in the time interval $[0,6]$ is $\Delta t=0.01$, and the space domain of QTT method is employed on $[-\frac{7}{2},\frac{7}{2}]^d$, with the degree of freedom $2^5$ in each direction. We implemented this example in $10$ and $20$-dimensional case. The number of particles in PF is chosen to maintain comparable relative mean squared error (RMSE) with TT approaches under the same setting (see 20D example).
%\begin{comment}
\begin{figure}[htbp]  
    \centering \label{fig:cubic_10D}
    \subfigure[1st dimension]{\includegraphics[width=0.32\linewidth]{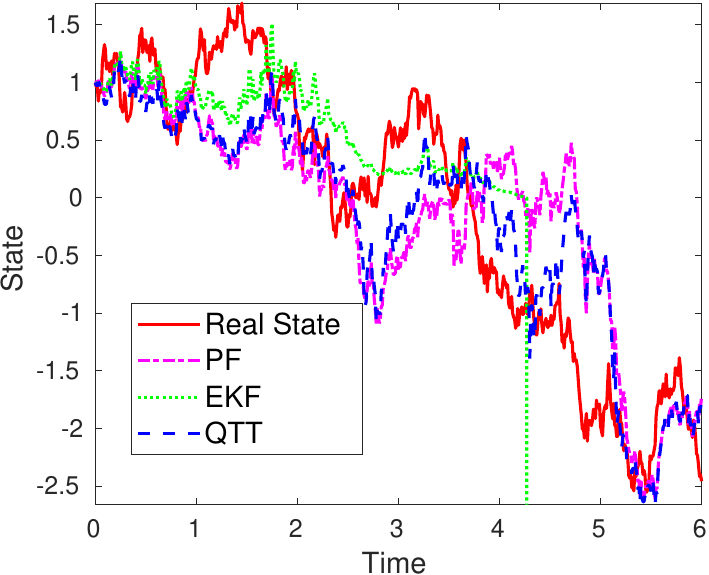}}
    \subfigure[2nd dimension]{\includegraphics[width=0.32\linewidth]{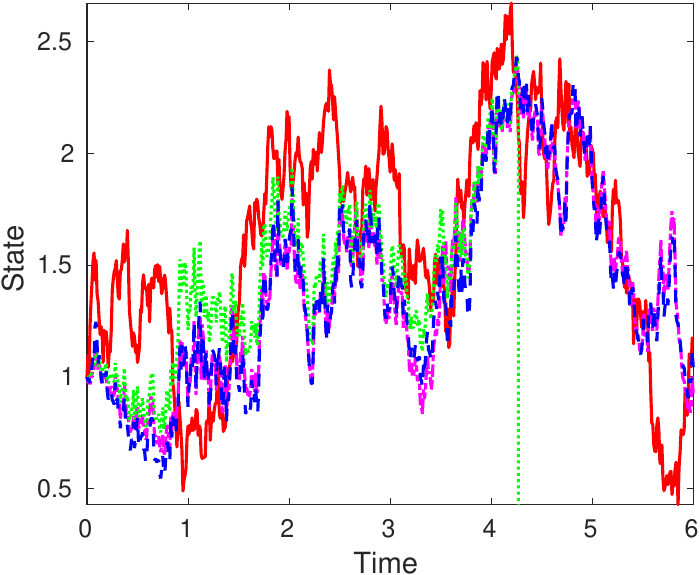}}
    \subfigure[3rd dimension]{\includegraphics[width=0.32\linewidth]{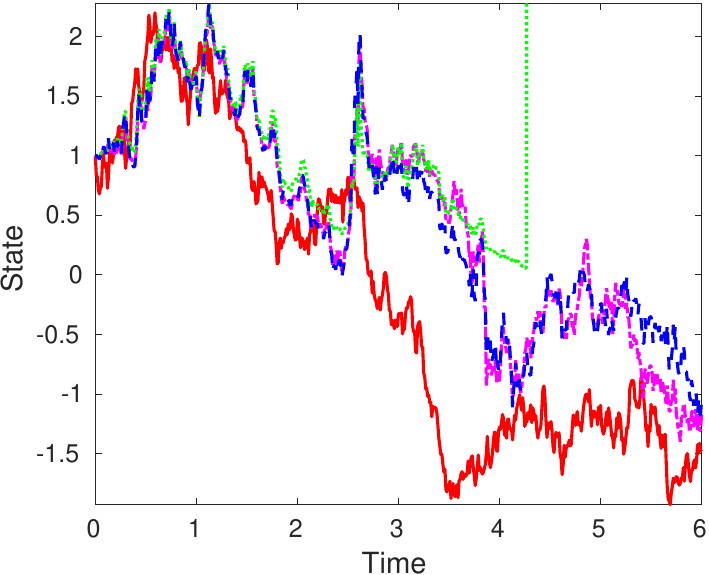}}
    \subfigure[4th dimension]{\includegraphics[width=0.32\linewidth]{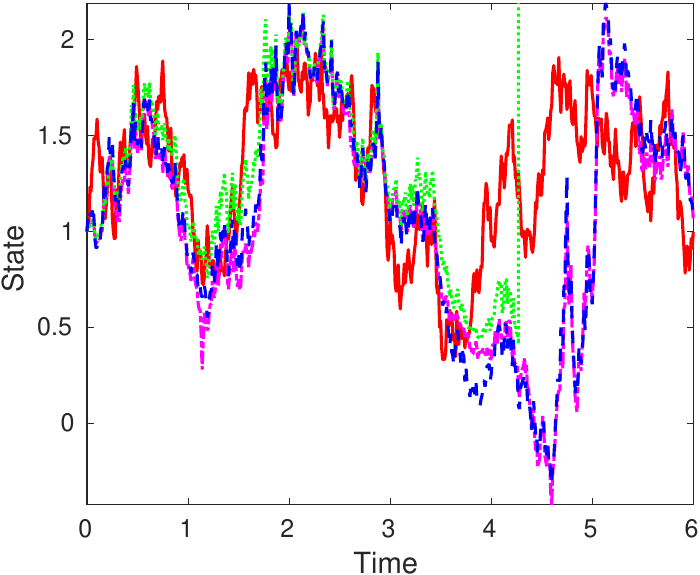}}
    \subfigure[5th dimension]{\includegraphics[width=0.32\linewidth]{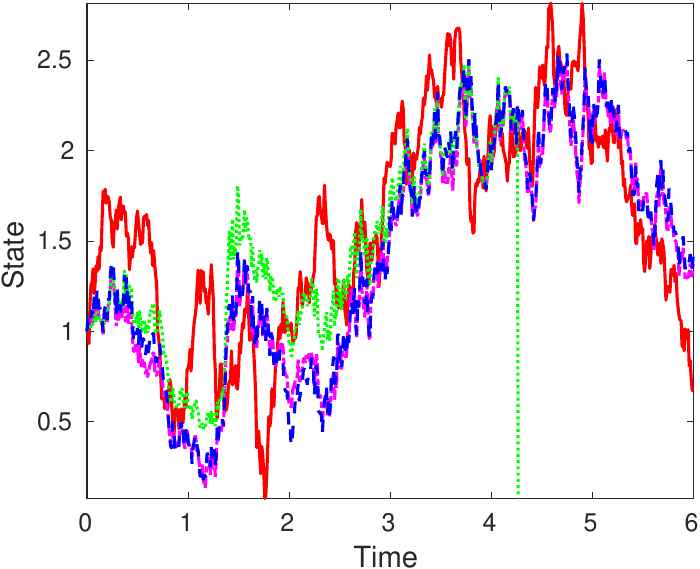}}
    \subfigure[6th dimension]{\includegraphics[width=0.32\linewidth]{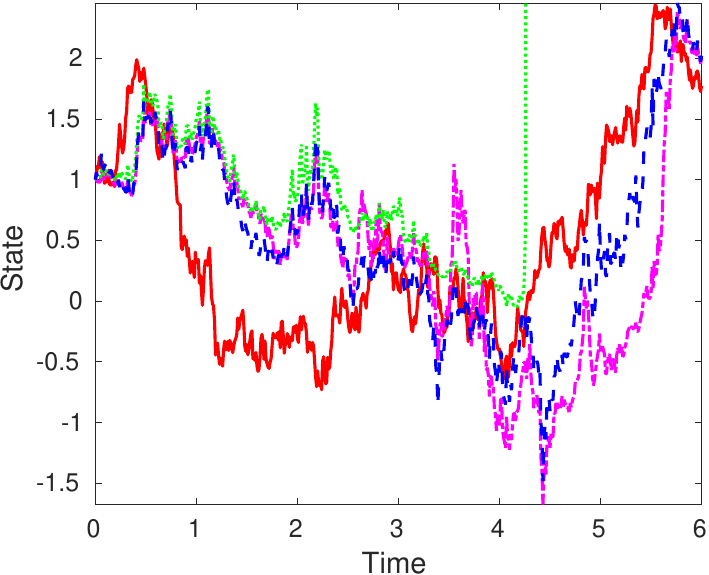}}
    \caption{Tracking results in 1-6 dimensions of QTT, PF (with 20000 particles), and EKF for the $10$D model \cref{eq:cubicsensor}.}
\end{figure}
%\end{comment}

In $10$D case, we employ $20000$ particles in PF and repeated the experiments $N_{path}=100$ times and observed that the RMSE of QTT and PF were $0.2757$ and $0.3015$, respectively. To achieve this level of accuracy, QTT required the average CPU time of $167.989$s, compared to $308.09$s for PF. Notably, the linearized filter, EKF, exhibited instability during tracking in most of the trials ($97$ out of $100$). 

\cref{fig:cubic_10D} illustrates the tracking results of $10$D dimensional case in a single experiment. In this trial, QTT, PF completed the task in $175.17$s and $341.37$s, respectively, both successfully track the real state. In contrast, the EKF diverged early ($t=4.2s$) and significantly deviated from the real trajectories. 

\begin{figure}[htbp]  
    \centering \label{fig:cubic_20D}
    \subfigure[1st dimension]{\includegraphics[width=0.32\linewidth]{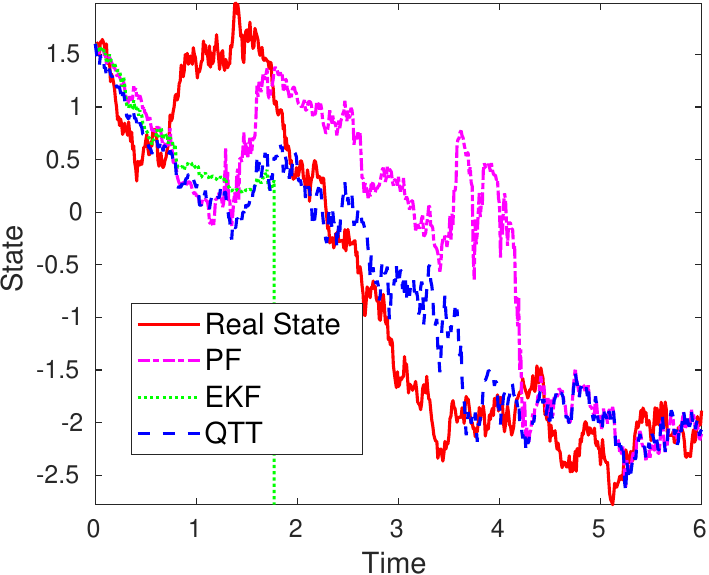}}
    \subfigure[2nd dimension]{\includegraphics[width=0.32\linewidth]{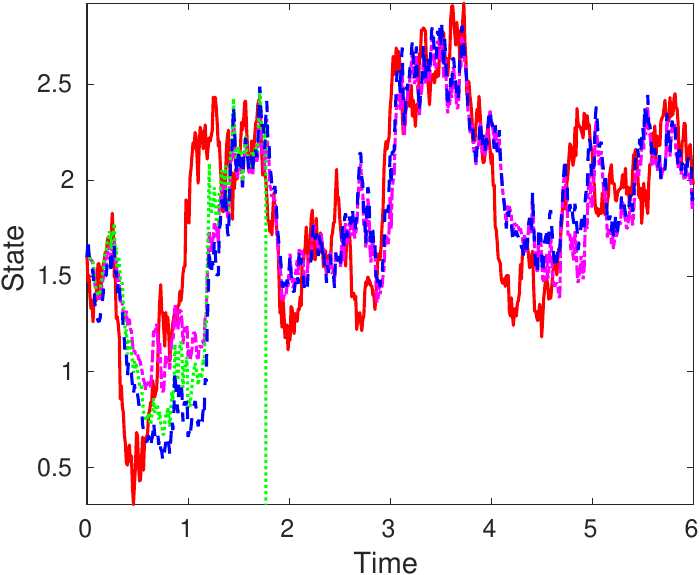}}
    \subfigure[3rd dimension]{\includegraphics[width=0.32\linewidth]{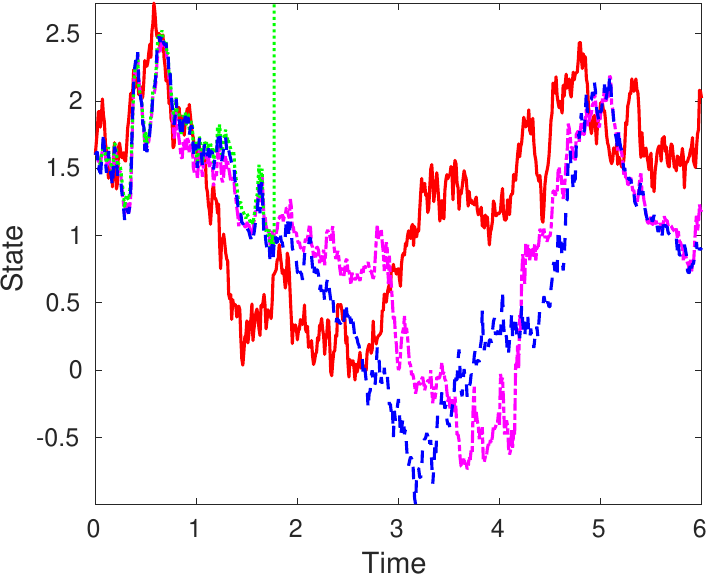}}
    \subfigure[4th dimension]{\includegraphics[width=0.32\linewidth]{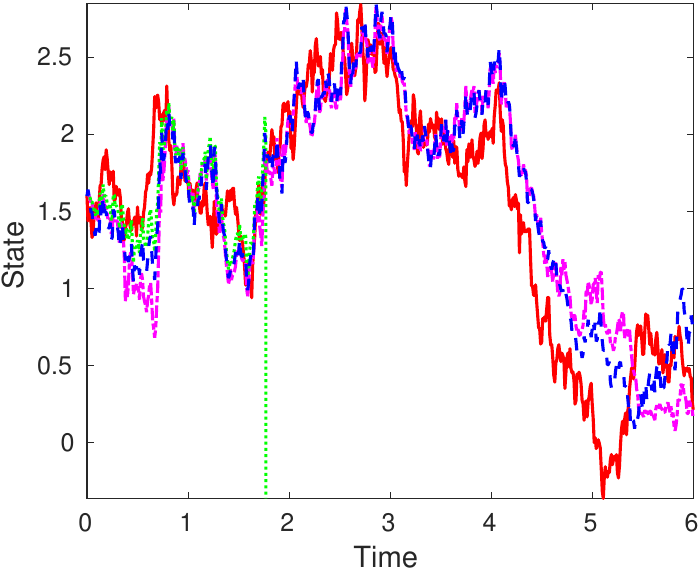}}
    \subfigure[5th dimension]{\includegraphics[width=0.32\linewidth]{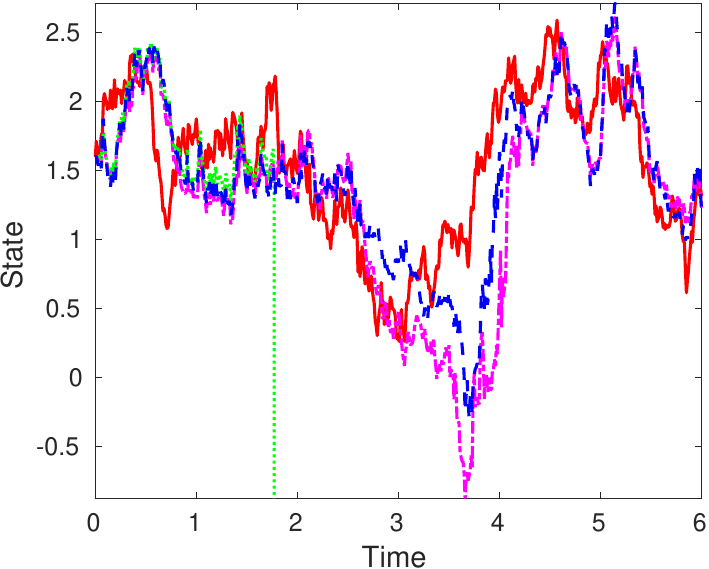}}
    \subfigure[6th dimension]{\includegraphics[width=0.32\linewidth]{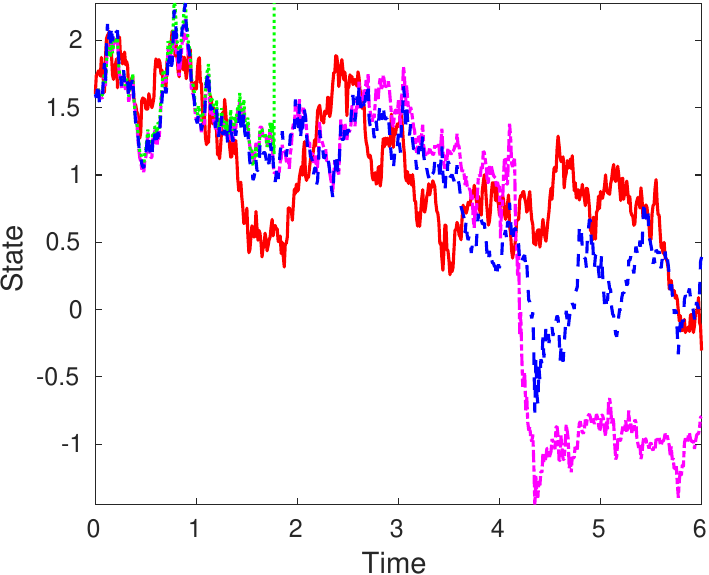}}
    \caption{Tracking results in 1-6 dimensions of QTT, PF (with 50000 particles), and EKF for the $20$D cubic sensor model \cref{eq:cubicsensor}.}
\end{figure}

In $20$D case, the initial state $x_0$ is set as $1.6\times{\bf 1}\in \mathbb{R}^{20}$ to enhance the signal-to-noise rate (SNR).  The tracking results of a single trial in $20$D case is depicted in \cref{fig:cubic_20D}. In this trial, the EKF exhibits numerical instability near $t=1.8s$. Furthermore, we repeat $100$ experiments under the same configuration with different paths of state and observation. In \cref{tab:cubci_20d}, we present the RMSE and computational cost of 600 steps of QTT and PF, where the computational cost of PF is $2.5\times$ of QTT for comparable RMSE, similar to $10$D model. To be noted, the EKF is unstable and deviates from the real trajectories in every trial out of the $100$ experiments.

\begin{table}[htbp]
\caption{Average RMSE and CPU time (online) of QTT and PF on the $20$D cubic sensor model}\label{tab:cubci_20d}
\begin{center}
  \begin{tabular}{ccc} \toprule
   Method & RMSE &  online time (sec) \\ \midrule
    QTT & 0.3184 & 532.15  \\
    PF (10000 particles) & 1.3622 & 241.66\\ 
    PF (30000 particles) & 0.9649 & 729.50 \\ 
    PF (50000 particles) & 0.4064 & 1287.16\\ \bottomrule
  \end{tabular}
\end{center}
\end{table}
%\begin{comment}
%\end{comment}
\subsection{Multi-mode problem}
Multi-mode conditional density functions may appear in some problems, e.g. the observe function $h(x_t)$ is an even function which is unable to recognize the sign of real trajectory. In such cases, EKF is more likely to fail with the assumption of a Gaussian distribution, while QTT can give a reasonable estimate by recovering the density function. 

Consider the multi-mode sensor problem modeled by :
\begin{equation}\label{bimodal}
    \begin{aligned}
\begin{cases}
&d x_1=(0.6\sin(x_1)+0.2\sin(x_2))dt+dv_1\\
&d x_2=(0.6\sin(x_2)+0.2\sin(x_3))dt+dv_2\\
&d x_3=0.6\sin(x_3)dt+dv_3\\
&dy_i=(x_i^2+\cos(x_i))+dt+dw_i,\quad i=1,2,3.\\
\end{cases}
\end{aligned}
\end{equation}

where $E(dv_tdv_t^T) = 0.2I_3dt, E(dw_tdw_t^T)=0.8 I_3$, and the initial state follows the Normal distribution $\mathcal{N}(0,\frac{1}{2\sqrt{5}}I_3)$. With the assigned drift function, observe function and the symmetric initial distribution, the density function $\sigma(x,t)$ is symmetric in $x$. The real states and observations are obtained in the time interval $[0, 10]$ with time step $\Delta t=0.01$. The space domain of solving FKE by QTT is $[-4,4]^3$, and the degree of freedom in each direction is $2^6$. The time step $\tau$ to solve local FKE is $\frac{\Delta t}{200}$. 

In the multi-mode problems, the expectation value may not be an appropriate statistic anymore, for instance, in problem \eqref{bimodal}, the ground truth of the conditional expectation will be $0$ due to symmetry. Therefore, we need to reconstruct the density function of the trajectories, to reflect potentially multiple high probability intervals.   In PF, one can use the normalized histogram of the simulated particles; in the proposed QTT appraoch, we use a heat map of the discretized density function, which is the full tensor recovered from the QTT approximation. %For a fair comparison of performance using the heat maps, the number of histogram bins is set to align with the degree of freedom in each direction in QTT method, keeping the colorbar scales nearly identical. 
%\begin{comment}
\begin{figure}
    \centering
    \subfigure[1st dimension(QTT)]{\includegraphics[width=0.32\linewidth]{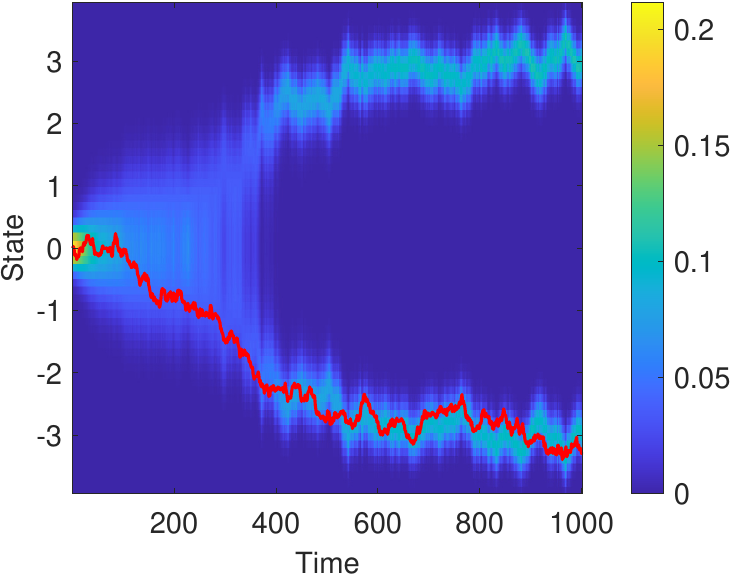}}
    \subfigure[2nd dimension(QTT)]{\includegraphics[width=0.32\linewidth]{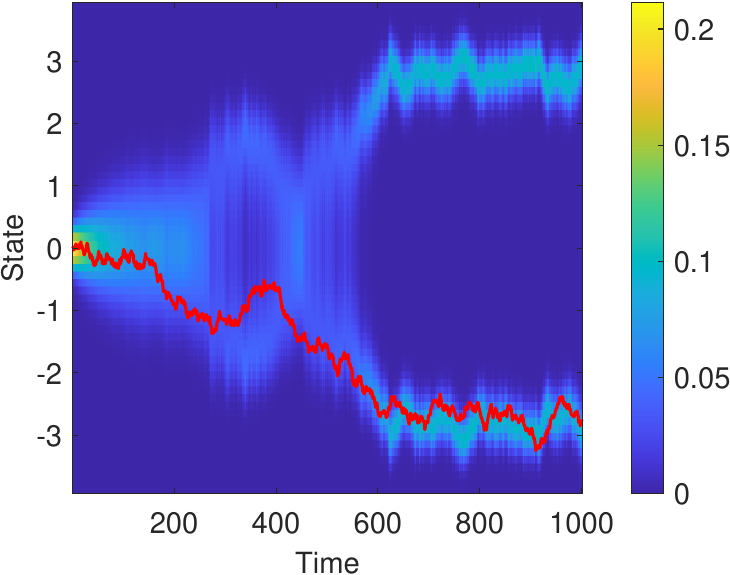}}
    \subfigure[3rd dimension(QTT)]{\includegraphics[width=0.32\linewidth]{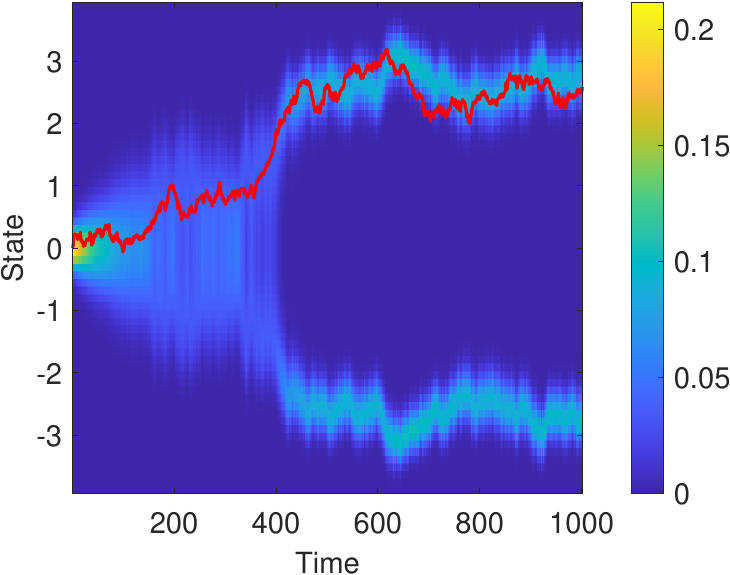}}
    \subfigure[1st dimension(PF)]{\includegraphics[width=0.32\linewidth]{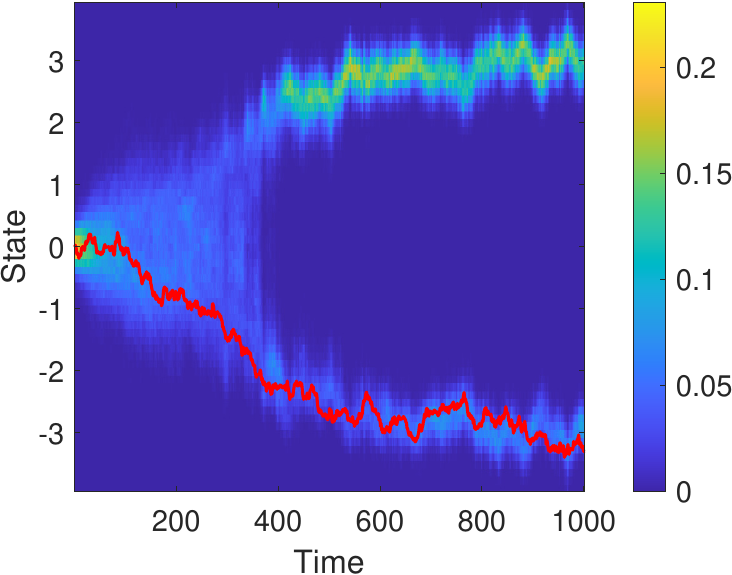}}
    \subfigure[2nd dimension(PF)]{\includegraphics[width=0.32\linewidth]{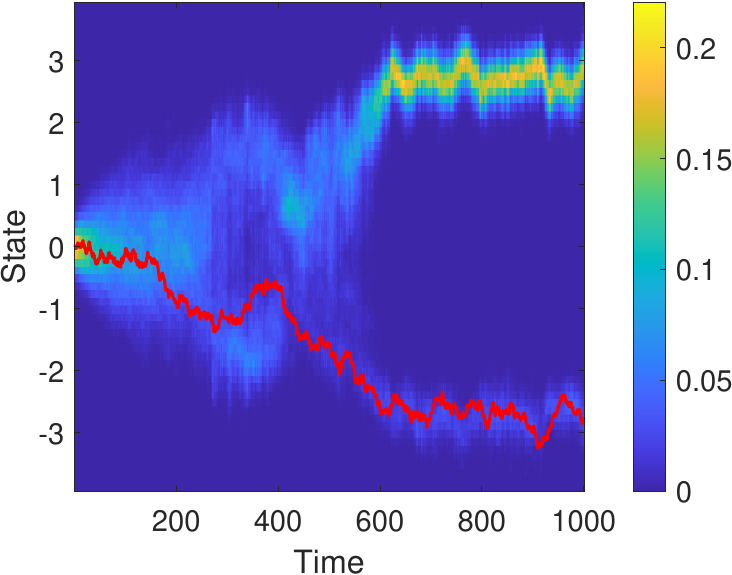}}	
    \subfigure[3rd dimension(PF)]{\includegraphics[width=0.32\linewidth]{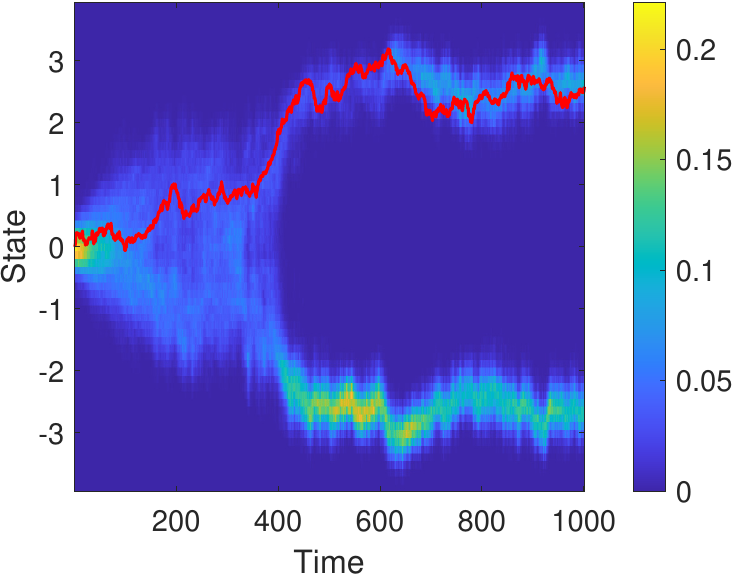}}
    \subfigure[1st dimension(EKF)]{\includegraphics[width=0.32\linewidth]{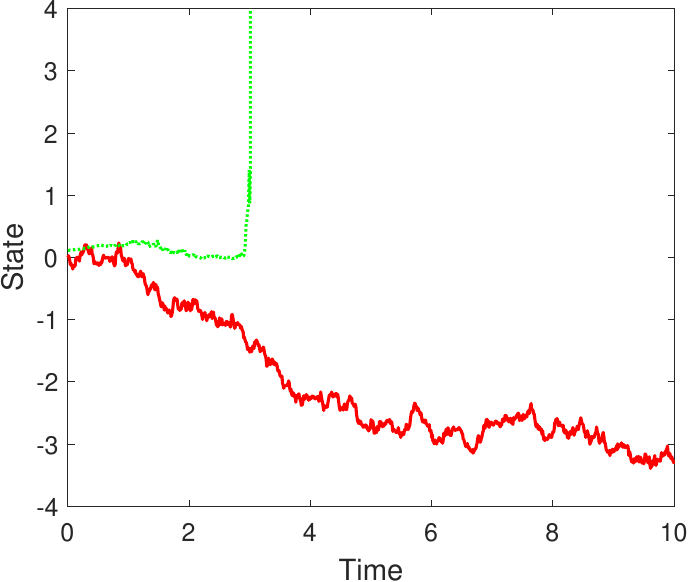}}
    \subfigure[2nd dimension(EKF)]{\includegraphics[width=0.32\linewidth]{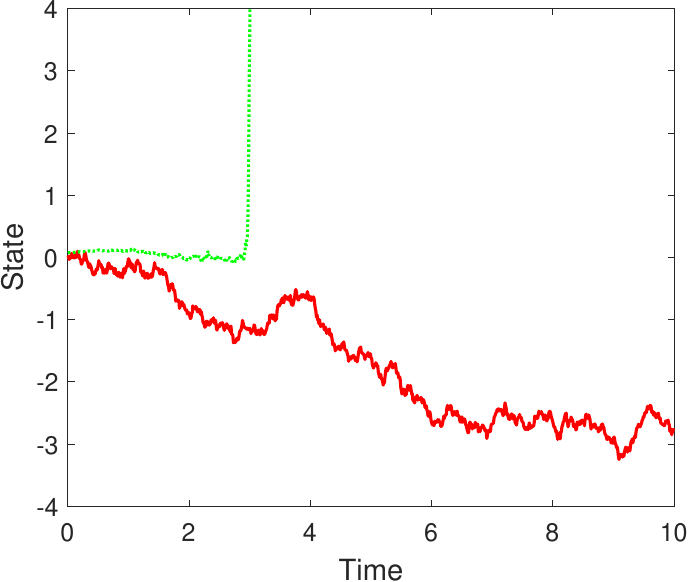}}	
    \subfigure[3rd dimension(EKF)]{\includegraphics[width=0.32\linewidth]{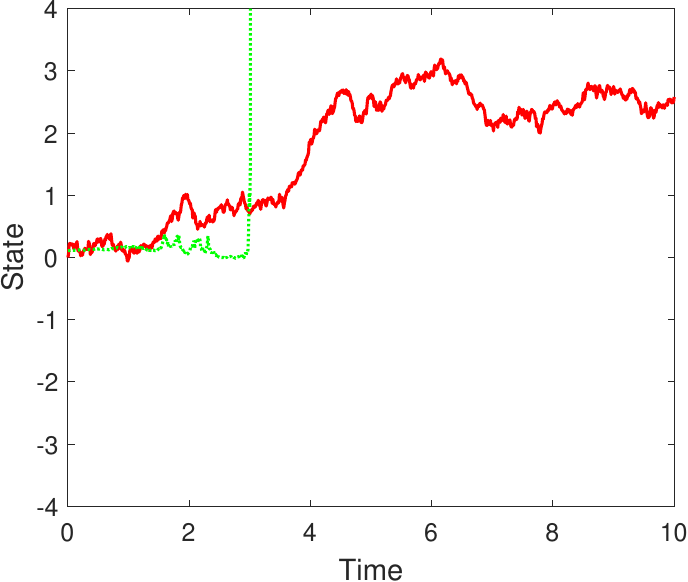}}
    \caption{Tracking results of QTT method, PF, and EKF for the model \cref{bimodal}. The red line is the ground truth state. In the first row (QTT), the background are the normalized density function; in the second row (PF), we plot the histogram of the particles (density normalization); in the third row (EKF), we plot the estimated state. } 
    \label{fig:bimodal}
\end{figure}
%\end{comment}

The tracking results of example \cref{bimodal} in a single trial is depicted in figure \ref{fig:bimodal}, the approximate density function obtained by QTT method reveals the bimodal and symmetric density function and closely follows the real state. PF method approximates the bimodal distribution but fails to achieve a symmetric distribution because of the randomness in sampling. EKF fails to track the state and encounters numerical instabilities.

\section{Conclusions}
\label{sec:conclusions}

In this paper, by using a weighted Sobolev space with power-type weights, we have presented the \textit{a priori} estimate of arbitrary order derivatives of the solution to the PR-DMZ equation in \cref{thm:prioriesti}. The embedding theorems then transform the estimate into an improved one, as stated in \cref{thm:improved_regularity_est}. The higher-order regularity estimate enhances the rationality of exploiting the low-rank structure and the convergence analysis of the QTT-based FD scheme. Equipped with the regularity results, we prove the convergence of this QTT method in \cref{thm:QTT_true_error_est}. More precisely, the error of the approximate solution arises from the prescribed tolerance in the TT operations and the mesh size in the FD scheme. From an implementation point of view, we show that the low storage and computation costs of the QTT approach to the DMZ equation can be ensured by assuming the functional polyadic structure of the drift and the observation terms. Finally, the numerical simulations demonstrate the efficiency and accuracy of our method compared with the existing algorithms, including the PF and EKF, see \cref{tab:cubci_20d}. Since our method reconstructs the conditional density function, it also gives an accurate estimate in multi-mode problems, see \cref{fig:bimodal}.

\appendix
%\section{Algorithm of PF and EKF}

\section{Functional Tensor Train decomposition}
\label{app:FTT decomposition}
The continuous analogue of discrete tensor train decomposition is derived in \cite{bigoni2016spectral}. Let $u\in L^2({\bf I})$, ${\bf I}=I_1\times I_2\times\cdots\times I_d$ be a bounded cube in $\mathbb{R}^d$. Denote that $\tilde{\bf{I}}:=I_2\times \cdot\cdot\cdot\times I_d$, then the operator
\begin{align*}
    T: L^2(\tilde{\bf{I}})&\longrightarrow L^2(I_1)\\
    g(x_2,...,x_d) &\longmapsto \int_{\tilde{\bf{I}}} u(x_1,x_2,...,x_d)g(x_2,...,x_d)dx_2\cdots dx_d
\end{align*}
is linear, bounded and compact. The Hilbert adjoint operator of $T$ is 
\begin{align*}
    T^*: L^2(I_1)&\longrightarrow L^2(\tilde{\bf{I}})\\
    f(x_1) &\longmapsto \int_{I_1} u(x_1,x_2,...,x_d)f(x_1)dx_1.
\end{align*}
Then $TT^*$ is a compact self-adjoint operator. By the theory of compact operator, the eigenvalues $\{\lambda_{1,\widetilde\alpha_1}\}_{\widetilde\alpha_1=1}^\infty$ of $TT^*$ are at most countable and the corresponding eigenfunctions $\{\varphi_1(x_1;\widetilde\alpha_1)\}_{\widetilde\alpha_1=1}^{\infty}$ form an orthonormal basis of $L^2(I_1)$. Similarly, the eigenfunctions $\{\phi_1({\widetilde\alpha_1};x_2,...,x_d)\}_{\widetilde\alpha_1=1}^{\infty}$ of $T^*T$ also form an orthonormal basis of $L^2(\tilde{\bf{I}})$. Consequently the function $u$ can be expressed in Schmidt decomposition
\begin{align}\label{Schmidt_decomposition}
    u=\sum_{\widetilde\alpha_1=1}^\infty \sqrt{\lambda_{1,\widetilde\alpha_1}}\varphi_1(x_1;\widetilde\alpha_1)\phi_1({\widetilde\alpha_1};x_2,...,x_d),
\end{align}
where the convergence is in $L^2({\bf I})$. Apply the same decomposition on $$\sqrt{\lambda_{1,\widetilde\alpha_1}}\phi_1({\widetilde\alpha_1};x_2,...,x_d),$$ the decomposition \cref{Schmidt_decomposition} is substituted by 
\begin{align*}
    u=\sum_{\widetilde\alpha_1=1}^\infty \sum_{\widetilde\alpha_2=1}^\infty \sqrt{\lambda_{2,\widetilde\alpha_2}}\varphi_1(x_1;\widetilde\alpha_1)\varphi_2(\widetilde\alpha_1;x_2;\widetilde\alpha_2)\phi({\widetilde\alpha_1};x_2,...,x_d).
\end{align*}
By recursion, the \textit{functional Tensor Train decomposition} follows 
\begin{align*}\label{FTT_decomposition}
    u = \sum_{\widetilde\alpha_1,...,\widetilde\alpha_{d-1}=1}^\infty \varphi_1(\widetilde\alpha_0;x_1;\widetilde\alpha_1) \varphi_2(\widetilde\alpha_1;x_2;\widetilde\alpha_2)\cdots\varphi_d(\widetilde\alpha_{d-1};x_d;\widetilde\alpha_d),
\end{align*}
where $\widetilde\alpha_0,\widetilde\alpha_d=1$. 

\section{Algorithm of PF}
The algorithm for particle filter is presented in \cref{alg:PF}.
\begin{algorithm}
\caption{Particle Filter Algorithm}
\begin{algorithmic}[1]\label{alg:PF}
    \Input Initial state $x_0$, initial distribution $\sigma_0(x)$, process noise covariance ${\bf Q}$, observe noise ${\bf S}=s{\bf Id}$, particles number $N$, drift function $f$, observe function $h$, observation $\{y_{t_{j}} \}_{j=0}^{N_T}$.
    \Output Estimated state trajectory $\{\hat x_{t_j}\}_{j=0}^{N_T}$.
    \State $\hat x_{t_0}\leftarrow x_0$.
    \State Draw initial particles $x_0^{(i)},i=1,...,N$ from initial distribution $\sigma_0(x)$.
    \State Initialize weights: $w_{0}^{(i)}=\frac{1}{N}$.
    \For {$j=1\rightarrow N_T$}
    \For {$i=1\rightarrow N$}
    \State {Prediction:} $\widetilde x_{t_j}^{(i)}=f(x_{t_{j-1}}^{(i)})\Delta t+\mathrm{d}v_{t_j}^{(i)}$, where $\mathrm{d}v_{j}^{(i)}\sim \mathcal{N}(0,{\bf Q}\Delta t)$.
    \State {Weight update: }Calculate likelihood: $$w_{j}^{(i)}=\exp\Big(-\frac {| y_{t_{j}}^{(i)}-y_{t_{j-1}}^{(i)}-h(\widetilde x_{t_j}^{(i)})\Delta t  |^2} {2s\Delta t}\Big).$$
    \EndFor
    \State {Normalize weights:} $w_{j}^{(i)}\leftarrow {w_{j}^{(i)}}/({\sum_k w_{j}^{(k)}}) $.
    \State \{Resampling\} 
    \State Construct cumulative density function: $c_i=\sum_{m=1}^jw_{j}^{(m)}, j=1,...,N$.
    \For {$i=1\rightarrow N$}
    \State Sample $r\sim\mathcal{U}(0,1)$.
    \State Find $k_i=\min\{j:c_j\ge r\}$.
    \State Set $x_{t_{j}}^{(i)}=\widetilde x_{t_j}^{(k_i)}$.
    \State Reset weights: $w_j^{(i)}=\frac{1}{N}$.
    \EndFor
    \State State estimate: $\hat x_{t_j}=\frac{1}{N}\sum_{i=1}^Nx_{t_j}^{(i)}$
    \EndFor 
\end{algorithmic}
\end{algorithm}

\section{Algorithm of EKF}
The algorithm for extended Kalman filter is presented in \cref{alg:EKF}.
\begin{algorithm}
\caption{Extended Kalman Filter Algorithm}
\begin{algorithmic}[1]\label{alg:EKF}
    \Input Initial state $x_0$, initial distribution $\mathcal{N}(x_0,{\bf P}_0)$ process noise covariance ${\bf Q}$, observe noise ${\bf S}$, particles number $N$, drift function $f$, observe function $h$, observation $\{y_{t_{j}} \}_{j=0}^{N_T}$
    \Output Estimated state trajectory $\{\hat x_{t_j}\}_{j=0}^{N_T}$.
    \State $\hat x_0\leftarrow x_{t_0}$.
    \For {$j=1\rightarrow N_T$}
    \State Compute Jacobian matrices: 
    $ {\bf F}_j=\frac{\partial f}{\partial x}(\hat x_{j-1}),\,{\bf H}_j = \frac{\partial h}{\partial x}(\hat x_{j-1})$.
    \State $\Delta {\bf P}_j=( {\bf F}_j{\bf P}_{j-1} + {\bf P}_{j-1}{\bf F}_j^T - {\bf P}_j({\bf H}_j^T{\bf S}^{-1}{\bf H}_j){\bf P}_j + {\bf Q} )\Delta t$.
    \State $\Delta\hat x_{t_j}=f(\hat x_{t_{j-1}})\Delta t$.
    \State ${\bf P}_j = {\bf P}_{j-1}+\Delta {\bf P}_{j}.$
    \State ${\bf K}_j={\bf P}_j{\bf H}_j^T{\bf S}^{-1}.$
    \State $\hat x_{t_j}=\hat x_{t_{j-1}}+\Delta\hat x_{t_j}$.
    \EndFor
\end{algorithmic}
\end{algorithm}

\newpage
\printbibliography
%\bibliographystyle{plain}   % or abbrv, siam, etc.
%\bibliography{references}  
\end{document}